\documentclass[10pt]{amsart}
\usepackage{amssymb,latexsym,amscd}
\input epsf
\makeatletter
\let\atopwithdelims\@@atopwithdelims
\let\over\@@over
\newtheorem{theorem}{Theorem}[section]
\newtheorem{proposition}[theorem]{Proposition}
\newtheorem{corollary}[theorem]{Corollary}

\newtheorem{remark}[theorem]{Remark}
\newtheorem{lemma}[theorem]{Lemma}

\newtheorem{question}[theorem]{Question}
\newcounter{enumit}

\newcommand{\reals}{{\mathbb R}}
\newcommand{\dimens}{\mathrm{dim}}
\newcommand{\Baues}{\omega}

\newcommand{\F}{{\mathcal F}}
\newcommand{\N}{{\mathcal N}}
\newcommand{\AAA}{{\mathcal A}}
\newcommand{\LL}{{\mathcal L}}
\newcommand{\B}{{\mathbb B}}
\newcommand{\SU}{{\mathbb S}}
\renewcommand{\to}{\rightarrow}
\newcommand{\toto}{\longrightarrow}

\newcommand{\sign} {\textrm{sign}}
\newcommand{\conv} {\mathrm{conv}}
\begin{document}
\title{Fiber polytopes for the projections between
cyclic polytopes}

\author{Christos~A.~Athanasiadis}
\address{\hskip-\parindent Christos~A.~Athanasiadis\\
Department of Mathematics\\
University of Pennsylvania\\
Philadelphia, PA 19104, USA}
\email{athana@math.upenn.edu}

\author{Jes\'us~ A.~De~Loera}
\address{\hskip-\parindent Jes\'us~A.~De Loera\\
The Geometry Center and School of Mathematics\\
University of Minnesota\\
Minneapolis, MN 55455, USA}
\email{deloera@geom.umn.edu}

\author{Victor Reiner}
\address{\hskip-\parindent Victor Reiner\\
School of Mathematics\\
University of Minnesota\\
Minneapolis, MN 55455, USA}
\email{reiner@math.umn.edu}

\author{Francisco Santos}
\address{\hskip-\parindent Francisco Santos\\
Departamento de Matem\'aticas, Estad\'{\i}stica y Computaci\'on \\
Universidad de Cantabria\\
Santander, E-39071, Spain}
\email{santos@matesco.unican.es}

\thanks{The first author was supported by a Mathematical Sciences
Research Institute postdoctoral fellowship and a University of Pennsylvania 
Hans Rademacher Instructorship. Research of the 
second and fourth authors were partially supported by the Geometry
Center (NSF grant DMS-8920161). The third author was supported
by a University of Minnesota McKnight-Land Grant Fellowship and Sloan 
Foundation Fellowship.}
\begin{abstract}
The cyclic polytope $C(n,d)$ is the convex hull of any $n$ points on
the moment curve $\{(t,t^2,\ldots,t^d):t \in \reals\}$ in $\reals^d$.
For $d' >d$, we consider the fiber polytope (in the sense of Billera and
Sturmfels) associated to the natural projection of cyclic polytopes
$\pi: C(n,d') \rightarrow C(n,d)$ which ``forgets" the last $d'-d$ coordinates.
It is known that this fiber polytope 
has face lattice indexed by the coherent polytopal subdivisions of
$C(n,d)$ which are induced by the map $\pi$. Our main result characterizes 
the triples $(n,d,d')$ for which
the fiber polytope is canonical in either of the following two senses: 
\begin{enumerate}
\item[$\bullet$] \textit{all} polytopal subdivisions induced by $\pi$ are
coherent,
\item[$\bullet$] the structure of the fiber polytope
does not depend upon the choice of points on the moment curve.
\end{enumerate}
We also discuss a new instance with a positive answer to the
Generalized Baues Problem, namely that of a projection $\pi:P\to Q$
where $Q$ has only regular subdivisions and $P$ has two more vertices
than its dimension.
\end{abstract}

\maketitle

\section{ Introduction}
\label{intro}

The cyclic $d$-polytope with $n$ vertices is the convex hull of any $n$ points
on the moment curve $\{(t,t^2,\ldots,t^d):t \in \reals\}$ in $\reals^d$.
Historically, the cyclic polytopes played an important role 
in polytope theory because they provide the upper bound for the number of
faces of a $d$-polytope with $n$ vertices
\cite[Chapter 8]{Ziegler},\cite[\S 4.7]{Grunbaum}.
Although the cyclic polytope itself depends upon the
choice of these $n$ points, much of its combinatorics,  
such as the structure of its lattice of faces or its set of triangulations,
is well-known to be independent of this choice 
(see \cite{Ziegler,EdRe,Rambau1}).
For this reason, we will often abuse notation and refer to the
cyclic $d$-polytope with $n$ vertices as $C(n,d)$, making reference
to the choice of points only when necessary.

The cyclic polytopes come equipped with a natural family of 
maps between them: fixing a pair of dimensions $d' > d$, the map 
$\pi:\reals^{d'}\rightarrow \reals^{d}$ which forgets the last
$d'-d$ coordinates restricts to a surjection 
$\pi: C(n,d')\rightarrow C(n,d)$.  Here we are implicitly
assuming that if the points on the moment curve in $\reals^{d'}$
chosen to define $C(n,d')$ have first coordinates $t_1 < \cdots < t_n$,
then the same is true for the points in $\reals^d$ chosen to define
$C(n,d)$.

Our starting point is that these maps $\pi: C(n,d')\rightarrow C(n,d)$
provide interesting and natural examples for 
Billera and Sturmfels' theory of \textit{fiber polytopes} \cite{BS1}.
Given an affine surjection of polytopes $\pi: P \rightarrow Q$,
the \textit{fiber polytope} $\Sigma(P \overset{\pi}{\rightarrow} Q)$ is a
polytope of dimension $\dimens(P) - \dimens(Q)$ which is (in a
well-defined sense; see \cite{BS1}) the ``average" fiber of the
map $\pi$.  The face poset of $\Sigma(P \overset{\pi}{\rightarrow} Q)$
has a beautiful combinatorial-geometric interpretation: it is
the refinement ordering on the set of all polytopal subdivisions of $Q$ which
are induced by the projection $\pi$ from $P$ in a certain combinatorial
sense, and also \textit{$\pi$-coherent} in a geometric sense -
see Section \ref{background}. This face poset sits inside the larger 
\textit{Baues poset} $\Baues(P \overset{\pi}{\rightarrow} Q)$, which is the 
refinement ordering on all polytopal subdivisions of $Q$ which are
induced by $\pi$, or {\it $\pi$-induced}. 

For the case of cyclic polytopes, the Baues poset 
$\Baues(C(n,d') \overset{\pi}{\rightarrow} C(n,d))$
of all $\pi$-induced subdivisions 
does not depend on the choice of points along the moment curve. On the
other hand, a $\pi$-induced subdivision may be $\pi$-coherent or not
depending on the choice of points.
The main question addressed by this paper is:  ``How canonical
is the fiber polytope  $\Sigma(C(n,d')\overset{\pi}{\rightarrow} C(n,d))$,
i.e. to what extent does its combinatorial structure vary with the choice
of points on the moment curve?". There are at least two ways in which
$\Sigma(C(n,d')\overset{\pi}{\rightarrow} C(n,d))$ can be canonical:

\begin{enumerate}
\item[$\bullet$] If all $\pi$-induced subdivisions of $C(n,d)$ are
$\pi$-coherent (for a certain choice of points, and hence for all
by Lemma \ref{oneall}) then the face lattice of
$\Sigma(C(n,d')\overset{\pi}{\rightarrow} C(n,d))$ coincides with the
Baues poset $\Baues(C(n,d') \overset{\pi}{\rightarrow} C(n,d))$.

\item[$\bullet$]
Even if there exist $\pi$-induced subdivisions
of $Q$ which are not $\pi$-coherent, it is possible that the
identity of the $\pi$-coherent subdivisions (and, in particular, the face
lattice of the fiber polytope) is independent of the choice of
points.
\end{enumerate}

Our main result characterizes exactly for which values of $n$, $d$ and
$d'$ each of these two situations occurs.

\begin{theorem}
\label{main-theorem}
Consider the map $\pi: C(n,d') \rightarrow C(n,d)$.
\begin{enumerate}
\item[(1)]
If $d=1$ then the set of $\pi$-coherent polytopal subdivisions of $C(n,1)$,
and hence the face lattice of the fiber polytope
$\Sigma(C(n,d')\overset{\pi}{\rightarrow} C(n,1))$,
is independent of the choice of points on the moment curve.
In fact, the face lattice of $\Sigma(C(n,d')\overset{\pi}{\rightarrow} C(n,1))$
coincides with that of the cyclic $(d'-1)$-zonotope having
$n-2$ zones.  Furthermore, all $\pi$-induced polytopal
subdivisions of $C(n,1)$ are $\pi$-coherent if and only if
$d'=n-1$ or $d' = 2$.
\item[(2)]
If $n-d'=1$ then:
\item[$\bullet$] If either $d \leq 2$ or $n-d \leq 3$
or $(n,d) \in \{(8,4),(8,3),(7,3)\}$ then all $\pi$-induced
subdivisions of $C(n,d)$ are $\pi$-coherent.
\item[$\bullet$] In all other cases with $n-d'=1$, there exists
a $\pi$-induced subdivision of $C(n,d)$ whose
$\pi$-coherence varies with the choice of points on
the moment curve and for every choice of points there is some
$\pi$-induced but not $\pi$-coherent subdivision.
\item[(3)]
If $d'-d=1$ then there are exactly two proper $\pi$-induced subdivisions,
both of them $\pi$-coherent in every choice of points.
\item[(4)]
If $n-d'\ge 2$, $d'-d\ge 2$ and $d\ge 2$ then
there exists a $\pi$-induced subdivision
of $C(n,d)$ whose $\pi$-coherence varies with the choice of points on
the moment curve and for every choice of points there is some
$\pi$-induced but not $\pi$-coherent subdivision.
\end{enumerate}
\end{theorem}

Part (1) is proved in Section \ref{monotone-path-polytopes}. 
The $\pi$-induced subdivisions in this case are the so-called \textit{
cellular strings} \cite{BKS} and the finest ones (the atoms in the
Baues poset) the \textit{monotone edge paths}. The fiber polytope in this case
is the so-called \textit{monotone path polytope} \cite{BS1, BKS}.
The \textit{cyclic zonotope} $Z(n,d)$, which appears in the statement,
is the Minkowski sum
of line segments in the directions of any $n$ points on the moment curve.
Like the cyclic polytope $C(n,d)$, its combinatorial
structure (face lattice) does not depend upon
the choice of points on the moment curve.

In the case of part (2) all subdivisions of $C(n,d)$ are
$\pi$-induced and the fiber polytope is the \textit{secondary polytope}
of $C(n,d)$, introduced by Gel'fand et al. \cite{GKZbook}. The same
authors \cite{GKZ1} and Lee \cite{Lee1,Lee2} proved that in the cases $d=1$ 
and $d=2$ or $n\le d+3$, respectively, all ($\pi$-induced) subdivisions are 
\textit{regular} for an arbitrary polytope $Q$.
Even more, it can be deduced from their work that the secondary
polytope of $C(n,d)$ is an $(n-2)$-cube for $d=1$, an $(n-3)$-dimensional
associahedron for $d=2$ and an $n$-gon for $n=d+3$.
We prove the rest of part (2) in Section \ref{secondary-polytopes}.

Part (3) is trivial and is included only for completeness. For any surjection
$\pi:P \rightarrow Q$ of a $(d+1)$-polytope $P$ onto a $d$-polytope $Q$
there are only two $\pi$-induced proper polytopal subdivisions, both
$\pi$-coherent:  the subdivisions of $Q$ induced by the ``upper" and 
``lower" faces of $P$ with respect to the projection $\pi$.

Part (4) is proved in Section \ref{counterexamples}.

Section \ref{Baues-section} deals with an instance of the Generalized Baues
Problem (GBP) posed in \cite{BKS}.
The GBP asks, in some sense, how close topologically
the Baues poset $\Baues(P \overset{\pi}{\rightarrow} Q)$ is to
the face poset of $\Sigma(P \overset{\pi}{\rightarrow} Q)$ 
inside it. The proper part of this face poset is the
face poset of the $(\dim(P) -\dim(Q)-1)$-dimensional
sphere which is the boundary of $\Sigma(P \overset{\pi}{\rightarrow} Q)$.
The GBP asks whether the proper part of the Baues poset
(suitably topologized \cite{Bjorner}) is homotopy equivalent to a 
$(\dim(P) -\dim(Q)-1)$-dimensional sphere.
This is known to be true when $\dim(Q)=1$ \cite{BKS} and when
$\dim(P)-\dim(Q) \leq 2$, but false in
general \cite{Rambau2,RamZie}. In previous work on cyclic polytopes 
(\cite{RaSa} and \cite{EdRaRe} for $d\le 3$)
it was shown to be true for 
$C(n,n-1) \overset{\pi}{\rightarrow} C(n,d)$.
We prove the following result, which in particular answers
the question positively for $C(n,n-2) \overset{\pi}{\rightarrow} C(n,2)$.

\begin{theorem}
\label{new-Baues-case}
Let $\pi: P \rightarrow Q$ have the property that
\begin{enumerate}
\item[$\bullet$] $P$ has $\dimens(P) + 2$ vertices and
\item[$\bullet$] the point configuration $\AAA$ which
is the image of the set of vertices of $P$ under $\pi$ has
only coherent subdivisions.
\end{enumerate}
Then the GBP has a positive answer for $\pi: P \rightarrow Q$.
\end{theorem}

One might be tempted to conjecture the following extension of Theorem
\ref{new-Baues-case}: the GBP has a positive answer if $\AAA$ has only
regular subdivisions, no matter what $P$ might be. However, one of the
counterexamples to the GBP given in \cite{RamZie} disproves this
extension. In that counterexample, $\AAA$ is planar and its 10 elements
are three copies of the vertices of a triangle, together with a point
inside.

\section{Background on fiber polytopes}
\label{background}

The fiber polytope $\Sigma(P\overset{\pi}{\rightarrow}Q)$,
introduced in \cite{BS1},
is a polytope naturally associated to any linear projection
of polytopes $\pi : P \rightarrow Q$.  An introduction to fiber polytopes
may also be found in \cite[Chapter 9]{Ziegler}.  In this section,
we review the definitions given in these two sources and discuss
some reformulations and further theory which we will need later.

Let $P$ be a $d'$-dimensional polytope with $n$ vertices in $\reals^{d'}$,
$Q$ a $d$-dimensional polytope in $\reals^d$ and 
$\pi: \reals^{d'} \toto \reals^d$ a linear map with $\pi(P)=Q$.
A \textit{polytopal subdivision} of $Q$ is a polytopal complex which
subdivides $Q$.  A polytopal subdivision of $Q$ is $\pi$-\textit{induced}
if
\begin{enumerate} 
\item[(i)] it is of the form $\{\pi(F):F \in \F\}$
for some specified collection $\F$ of faces of $P$ and 
\item[(ii)] $\pi(F) \subseteq \pi(F')$ implies $F = F' \cap \pi^{-1}(\pi(F))$,
thus in particular $F \subseteq F'$.
\end{enumerate}
It is possible that different collections $\F$ of faces of $P$ project
to the same subdivision $\{\pi(F):F \in \F\}$ of $Q$, so we distinguish
these subdivisions by labelling them with the family $\F$. Note that 
condition (ii) is superfluous for the family of projections $C(n,d')
\to C(n,d)$. We partially order the
$\pi$-induced subdivisions of $Q$ by $\F_1 \leq \F_2$ if and only if
$\bigcup \F_1 \subseteq \bigcup \F_2$. The resulting partially
ordered set is denoted by $\Baues(P {\overset{\pi}{\rightarrow}} Q)$
and called the \textit{Baues poset}.
The minimal elements in this poset are the \textit{tight} subdivisions,
i.e. those for which $F$ and $\pi(F)$ have the same dimension for all
$F$ in $\F$.

There are a number of ways to define $\pi$-coherent subdivisions of $Q$.
We start with the original definition from \cite{BS1}.
Choose a linear functional
$f \in (\reals^{d'})^*$. For each point $q$ in $Q$, the
fiber $\pi^{-1}(q)$ is a convex polytope which has a unique
face ${\overline F}_q$ on which the value of $f$ is minimized.
This face lies in the relative interior of a unique
face $F_q$ of $P$ and the collection of faces $\F=\{F_q\}_{q \in Q}$
projects under $\pi$ to a subdivision of $Q$.
Subdivisions of $Q$ which arise from
a functional $f$ in this fashion are called {\it $\pi$-coherent}.

It is worth mentioning here a slight variant of this description
of the $\pi$-coherent subdivision induced by $f$ (see also the paragraph
after the proof of Theorem 2.1 in \cite{BS2}). Note that
the inclusion $ker(\pi) \hookrightarrow \reals^{d'}$
induces a surjection $(\reals^{d'})^* \twoheadrightarrow ker(\pi)^*$.
A little thought shows that two functionals
$f,f' \in (\reals^{d'})^*$ having the same image
under this surjection will induce the same $\pi$-coherent
subdivision.  As a consequence, we may assume that the
functional $f$ lies in $ker(\pi)^*$.  From this point of view,
the following lemma should be clear.

\begin{lemma}
\label{coherence-by-normal-cones}
A face $F$ of $P$ belongs to the $\pi$-coherent subdivision of $Q$
induced by $f\in ker(\pi)^*$ if and only if its {\it normal cone}
$\N(F) \subseteq (\reals^{d'})^*$ has the property that its image
under the surjection $(\reals^{d'})^* \twoheadrightarrow
ker(\pi)^*$ contains $f$.
\end{lemma} 

In \cite[\S 9.1]{Ziegler}, Ziegler defines $\pi$-coherent subdivisions
in the following
equivalent fashion.  Having chosen the functional $f \in (\reals^{d'})^*$
as above, form the graph of the linear map $\hat{\pi}:P \toto \reals^{d+1}$
given by $p \mapsto (\pi(p), f(p))$. The image of this map is a polytope 
$\hat{Q}$ in $\reals^{d+1}$ which maps onto $Q$ under the projection 
$\reals^{d+1} \rightarrow \reals^d$ which forgets the last coordinate.
Therefore, the set of \textit{lower faces} of $\hat{Q}$
(those faces whose normal cone contains a vector with negative
last coordinate) form a polytopal subdivision of $Q$. We identify
this subdivision of $Q$ with the family of faces $\F = \{F\}$ in $P$ which
are the inverse images under $\hat{\pi}$ of the lower faces of $\hat{Q}$.
Under this identification, it is easy to check that the subdivision
of $Q$ is exactly the same as the $\pi$-coherent subdivision induced by
$f$, described earlier. Let $\Baues_{coh}(P \overset{\pi}{\rightarrow} Q)$
denote the induced subposet of $\Baues(P \overset{\pi}{\rightarrow} Q)$
consisting of all $\pi$-coherent subdivisions of $Q$.

\begin{theorem}\cite{BS1}
The poset $\Baues_{coh}(P \overset{\pi}{\rightarrow} Q)$ is the face lattice
of a $(d'-d)$-dimensional polytope, the {\it fiber polytope}
$\Sigma(P \overset{\pi}{\rightarrow} Q)$.
\end{theorem} 

It will be useful for us later to have a reformulation of
these definitions using \textit{affine
functionals}, \textit{Gale transforms} and \textit{secondary polytopes}.
For this purpose, given
our previous situation of a linear map of polytopes $\pi:P \to Q$,
define a map $\phi_P: \reals^n \to \reals^{d'+1}$ by the
$(d'+1) \times n$ matrix having the vertices $p_i$ of $P$
as its columns and an extra row on top consisting of all ones.
Let $q_i = \pi(p_i) \in Q$ and define similarly the map 
$\phi_Q: \reals^n \toto \reals^{d+1}$. Then $\pi$ extends to a map
$\pi:\reals^{d'+1} \toto \reals^{d+1}$ such that $\pi \circ \phi_P
=\phi_Q$. In other words, the following diagram commutes:
\begin{equation}
\label{diagram}
\begin{CD}
\reals^n   @>{\phi_P}>>  \reals^{d'+1} \\
&    \searrow^{\phi_Q }{}   & {\pi}{}{\Bigg\downarrow}\\
     & & \reals^{d+1} \\
\end{CD}
\end{equation}

Consider the map $\phi_Q: \reals^n \toto \reals^{d+1}$ as a
projection onto $Q$ of the $(n-1)$-simplex $\Delta^{n-1}$ whose vertices
are the standard basis vectors in $\reals^n$. Given a linear
functional $f \in (\reals^n)^*$, we can interpret the $\phi_Q$-coherent 
subdivision of $Q$ induced by $f$ in the following fashion, using Ziegler's 
description: Write $f({\bf x}) = \sum_i w_i x_i$ and lift the $i^{th}$
vertex in $Q$ (i.e. the image under $\phi_Q$
of the $i^{th}$ standard basis vector) into $\reals^{d+1}$ with
last coordinate $w_i$. Then take the convex hull of these points to form
a polytope $\hat{Q}$. The lower faces of $\hat{Q}$ form
the desired $\phi_Q$-coherent subdivision, which is sometimes
referred to as the \textit{regular subdivision} induced by the
heights $w_i$.

We wish to describe when two functionals $f, f'$ induce
the same regular subdivision.  As before, this will certainly
be true whenever they have the same image under the surjection 
$(\reals^n)^* \twoheadrightarrow ker(\phi_Q)^*$ and therefore
one may consider them as elements of $ker(\phi_Q)^*$.
Let $G_Q$ be any $(n-d-1) \times n$ matrix whose rows form a basis
for $ker(\phi_Q)$. The {\it Gale transform} $Q^*$ is defined
to be the vector configuration $q_1^*,\ldots,q_n^*$ given by the columns of
$G_Q$. Note that by construction, the row space $Row(G_Q)$
is identified with $ker(\phi_Q)$ and since there is
a canonical identification of the dual of the row space with the
column space, we have that $f$ is a vector in the column space 
$Col(G_Q)$, i.e. the space containing the Gale transform points
$q_1^*,\ldots,q_n^*$. The following lemma is a form of
oriented matroid duality (see \cite[Lemma 3.2]{BGS},\cite[Section 4]{BFS}).
\begin{lemma}
\label{OM-duality-lemma}
A subset $S \subseteq \{q_1,\ldots,q_n\}$
spans a subpolytope of $Q$ which appears in the
regular subdivision induced by $f$ if and only if $f$ lies
in the relative interior of the positive cone spanned by
$S^c:=\{q_i^*: q_i \not\in S\}$.
\end{lemma}

\noindent
Consequently, two functionals 
$f,f' \in ker(\phi_Q)^* = Col(G_Q)$ induce the same regular subdivision
if and only if they lie in the same face of the \textit{chamber complex} of 
$Q^*$, which is the common refinement of all positive cones spanned by
subsets of the Gale transform $Q^* = \{q^*_1,\ldots,q^*_n\}$.  The
chamber complex of $Q^*$ turns out to be the normal fan of the
\textit{secondary polytope}
$\Sigma(Q):=\Sigma(\Delta^{n-1}\overset{\phi_Q}{\to}Q)$
defined by Gel'fand, Kapranov and Zelevinsky \cite{GKZbook}.

It is also possible to determine in this picture when the regular subdivision
of $Q$ induced by $f \in (\reals^n)^*$ is also $\pi$-induced 
(see \cite[Theorem 2.4]{BS1}). The surjection
$\phi_P: \reals^n \twoheadrightarrow \reals^{d'+1}$ induces an inclusion
$\phi_P^*: (\reals^{d'+1})^* \hookrightarrow (\reals^{n})^*$.
It then follows immediately from Ziegler's description 
\cite[\S 9.1]{Ziegler}
that the $\pi$-induced subdivision of $Q$ induced by some functional
$f \in (\reals^{d'+1})^*$ is the same as the regular subdivision
of $Q$ induced by the functional 
$\phi_P^*(f)=f \circ \phi_P \in (\reals^{n})^*$.
In other words, a set of heights $(w_1,\ldots,w_n)$ induces
a regular subdivision of $Q$ which is also $\pi$-induced if and
only if $\sum_i w_i x_i \in im(\phi_P^*)$.  Since the vectors of
$im(\phi_P^*)$ are characterized by the fact that they vanish
on $ker(\phi_P)$ (``row-space is orthogonal to null-space"),
to check $\sum_i w_i x_i \in im(\phi_P^*)$ in
practice one only needs to verify that every affine dependence
$\sum_i c_i p_i= 0$ of the vertices of $P$ satisfies
$\sum_i c_i w_i = 0$.  We state this as a lemma for later use.

\begin{lemma}
\label{Paco-characterization}
For the projection $\pi: P \rightarrow Q$, a
regular subdivision of $Q$ is also $\pi$-coherent if and only if it
can be induced by a functional $f({\bf x}) = \sum_i w_i x_i \in (\reals^n)^*$
which vanishes on $ker(\phi_P)$ or, equivalently, which satisfies 
$\sum_i c_i w_i = 0$ for every affine dependence $\sum_i c_i p_i= 0$
of the vertices of $P$.
\end{lemma}

Furthermore, we
can identify the normal fan $\N(\Sigma(P \overset{\pi}{\to} Q))$
to the fiber polytope as a subspace intersected with the chamber
complex of $Q^*$. Composing the embedding
$(\reals^{d'+1})^* \hookrightarrow (\reals^n)^*$
with our earlier surjection
$(\reals^n)^* \twoheadrightarrow ker(\phi_Q)^*$
gives a map $\phi^*_{P,Q}$ whose image $im(\phi^*_{P,Q})$ is a subspace in
$ker(\phi_Q)^*=Col(G_Q)$.  The next proposition then follows
immediately from our previous discussion.
\begin{theorem}
\label{subspace-coherence}
A regular subdivision of $Q$ is $\pi$-coherent if and only if
the relative interior of its corresponding cone in the chamber complex of
$Q^*$ in $Col(G_Q)$ contains a vector in the subspace $im(\phi^*_{P,Q})$.
The normal fan $\N(\Sigma(P \overset{\pi}{\rightarrow} Q))$ is identified
with the cone complex obtained by intersecting the chamber complex
of $Q^*$ in $Col(G_Q)$ with the subspace $im(\phi^*_{P,Q})$.
\end{theorem}

Because the cyclic polytopes $C(n,d)$, $C(n+1,d+1)$ are related
by a single element lifting, it will later be necessary for
us to recall how the chamber complex and Gale transform behave
with respect to such liftings (see \cite[\S 3]{BGS}).  
Given a $d$-polytope $Q$ in $\reals^d$ with $n$ vertices
$q_1,\ldots,q_n$, we say that a $(d+1)$-polytope $\hat{Q}$ in $\reals^{d+1}$
is a \textit{single element lifting} of $Q$ if it has $n+1$ vertices
$\hat{q}_1,\ldots,\hat{q}_n,\hat{q}_{n+1}$ and there is surjection 
$f: \reals^{d+1} \rightarrow \reals^d$ satisfying $f(\hat{q}_{n+1})=0$
and $f(\hat{q}_i) = c_i q_i$ for $i \leq n$ and some positive
scalars $c_i \in \reals$.  For the case of $Q=C(n,d)$, $\hat{Q}=C(n+1,d+1)$,
if we assume that the parameters for the points on the moment curve are chosen
so that $t_1 < \cdots <t_n < t_{n+1}=0$, the map $f$ is the one
which ignores the first coordinate and reverses the signs of the rest
and the constant $c_i$ is $-t_i$.
When we have a single element lifting $\hat{Q}$ of $Q$, let
$\tau: \reals^{n+1} \rightarrow \reals^n$ be the map which
sends $e_{n+1}$ to $0$ and $e_i$ to $c_i e_i$ for $i \leq n$.
One can check that this definition makes the diagram
\begin{equation*}
\begin{CD}
   \reals^{n+1}   @>{\phi_{\hat{Q}}}>>  \reals^{d+2}  \\
    @V{\tau}VV                              @VVfV     \\
    \reals^n     @>{\phi_Q}>>              \reals^{d+1}    &  \\
\end{CD}
\end{equation*}
commute, where $\phi_Q, \phi_{\hat{Q}}$ were defined earlier and
the map $f$ has been extended  to $\reals^{d+2} \rightarrow \reals^{d+1}$
by mapping $e_1 \mapsto e_1$. One can check that under these hypotheses,
$\tau$ restricts to an isomorphism 
$\tau: ker(\phi_{\hat{Q}}) \rightarrow ker(\phi_Q)$ and hence
its dual $\tau^*$ gives an isomorphism 
$ker(\phi_Q)^* \rightarrow ker(\phi_{\hat{Q}})^*$ between the
spaces containing the Gale transforms  $Q^*, \hat{Q}^*$.

\begin{lemma}\cite[Lemma 3.4]{BGS}
\label{lifting-lemma}
When $\hat{Q}$ is a single-element lifting of $Q$,
one can choose the matrices $G_Q, G_{\hat{Q}}$ 
whose columns give the Gale transform points $Q^*, \hat{Q}^*$
in such a way that the isomorphism $\tau^*$ maps $q^*_i$ to $\hat{q}^*_i$
for $i \leq n$.
\end{lemma}
In other words, the Gale transform $\hat{Q}^*$ is a single-element 
extension of the Gale transform $Q^*$. 

We also wish to deal with the relation between the fiber polytopes for the
natural projections $\pi:C(n,d') \rightarrow C(n,d)$ and
$\hat{\pi}:C(n+1,d'+1) \rightarrow C(n+1,d+1)$. More generally,
given two single-element liftings $\hat{P}$ of $P$, with
map $f_P$, and $\hat{Q}$ of $Q$, with map $f_Q$,
we say that two linear surjections $\pi:P \rightarrow Q$
and $\hat{\pi}:\hat{P} \rightarrow \hat{Q}$ are \textit{compatible
with the liftings} if $f_Q \circ \hat\pi =\pi \circ f_P$.
Since $\phi_Q=\pi \circ \phi_P$ and
$\phi_{\hat Q}=\hat\pi \circ \phi_{\hat P}$, it easily follows that
$\tau_P=\tau_Q$ and the following diagram commutes:
\begin{equation}
\label{first-prism}
\begin{CD}
\reals^{n+1} @>{\phi_{\hat P}}>> \reals^{d'+2}  @>{\hat{\pi}}>> 
       \reals^{d+2}  @<{\phi_{\hat Q}}<< \reals^{n+1}\\
   @VV{\tau_P}V   @VV{f_P}V   @VV{f_Q}V   @VV{\tau_Q}V \\
\reals^{n} @>{\phi_{P}}>> \reals^{d'+1}  @>{\pi}>> 
       \reals^{d+1}  @<{\phi_{Q}}<< \reals^{n}\\
\end{CD} 
\end{equation}
One can easily check that this is the case for
$\pi:C(n,d') \rightarrow C(n,d)$ and
$\hat{\pi}:C(n+1,d'+1) \rightarrow C(n+1,d+1)$.

\begin{lemma}
\label{fiber-lifting-lemma}
Assume as above that $\hat{P}, \hat{Q}$ are single-element liftings of
$P,Q$ and that $\pi,\hat{\pi}$ are compatible projections. Then the
isomorphism $\tau^*:ker(\phi_Q) \rightarrow ker(\phi_{\hat Q})$
restricts to an isomorphism of the subspaces $im(\phi^*_{P,Q})$, 
$im(\phi^*_{\hat{P},\hat{Q}})$, which contain the normal fans of the
fiber polytopes $\Sigma(P \overset{\pi}{\rightarrow} Q)$,
$\Sigma(\hat{P} \overset{\hat{\pi}}{\rightarrow} \hat{Q})$.
\end{lemma}
\begin{proof}
As a preliminary step,
we note that the subspace $im(\phi^*_{P,Q})$ in $ker(\phi_Q)^*$
can be characterized slightly differently. Consider the short exact
sequence 
$$
\begin{CD}
0 \rightarrow ker(\phi_P) @>{i_{P,Q}}>> ker(\phi_Q) 
               @>{\phi_{P,Q}}>> \reals^{d'+1} \rightarrow 0,
\end{CD}
$$
where $\phi_{P,Q}$ is the composite 
$ker(\phi_Q) \hookrightarrow \reals^n \overset{\pi}{\twoheadrightarrow}
\reals^{d'+1}$ and $i_{P,Q}$ is the inclusion map.
Dualizing this says that $im(\phi^*_{P,Q}) = ker( i^*_{P,Q})$, which is
the characterization we will need.

The commutative diagram in (\ref{first-prism}) gives rise
to the following commutative square
\begin{equation*}
\begin{CD}
 ker(\phi_{\hat{P}}) @>{i_{\hat{P},\hat{Q}}}>>  ker(\phi_{\hat{Q}}) \\
 @V{\tau}VV   @VV{\tau}V \\
 ker(\phi_P)    @>{i_{P,Q}}>>ker(\phi_Q)\\ 
\end{CD}
\end{equation*}
Dualizing this square gives a square in which the kernels of the
horizontal maps can be added:
\begin{equation*}
\begin{CD}
ker(i^*_{\hat{P},\hat{Q}}) \hookrightarrow &
ker(\phi_{\hat{Q}})^* @>{i^*_{\hat{P},\hat{Q}}}>> ker(\phi_{\hat{P}})^*   \\
&   @A\tau^*AA      @AA\tau^*A  \\
ker(i^*_{P,Q})  \hookrightarrow   & ker(\phi_Q)^* 
@>{i^*_{P,Q}}>> ker(\phi_P)^*
\end{CD}
\end{equation*}
One then checks that the vertical map $\tau^*$ restricts to an isomorphism
$ker(i^*_{P,Q}) \rightarrow ker(i^*_{\hat{P},\hat{Q}})$. We combine this
with the fact that $ker(i^*_{P,Q}) = im(\phi^*_{P,Q})$ to get the
assertion.
\end{proof}

\section{The case $d=1$:  monotone path polytopes}
\label{monotone-path-polytopes}

In this section we restrict our attention to the natural  
projection $\pi: C(n,d') \to C(n,1)$ and  prove the assertions
of Theorem \ref{main-theorem} concerning the case $d=1$.
We recall and separate out these assertions in the following theorem
where, for ease of notation, we have replaced $d'$ by $d$. 

\begin{theorem}
\label{monopp}
For the natural projection $\pi: C(n,d) \to C(n,1)$,
the set of $\pi$-coherent polytopal subdivisions of $C(n,1)$,
and hence the face lattice of the fiber polytope
$\Sigma(C(n,d)\overset{\pi}{\rightarrow} C(n,1))$,
is independent of the choice of points on the moment curve.
In fact, the face lattice of 
$\Sigma(C(n,d)\overset{\pi}{\rightarrow} C(n,1))$
coincides with that of the cyclic $(d-1)$-zonotope having
$n-2$ zones.  Furthermore, all $\pi$-induced polytopal subdivisions of $C(n,1)$
are $\pi$-coherent if and only if $d=2$ or $d=n-1$.
\end{theorem}

In the case of $\pi:P \to Q$ with $\dim(Q)=1$, tight
$\pi$-coherent subdivisions $\F$ correspond to certain
\textit{monotone} edge paths on $P$.  The fiber polytope in this case
is called the \textit{monotone path polytope} (see \cite[\S 5]{BS1}
\cite[\S 9.2]{Ziegler} \cite{Ath} for examples). Before proving Theorem
\ref{monopp}, we recall the definitions of cyclic polytopes and
zonotopes, and describe explicitly the face lattice of the cyclic
zonotopes.

The {\it cyclic $d$-polytope with $n$ vertices}
$C(n, d)$ is the convex hull of the points
$$\{  v_i\} = \{(t_i, t_i^2,\ldots,t_i^d) \}_{i=1}^n $$
in ${\mathbb R}^d$, where $t_1 < t_2 < \cdots < t_n$.
Similarly, the \textit{cyclic $d$-zonotope with $n$ zones}  $Z(n,d)$ is the
$d$-zonotope generated by the vectors
$$\{ u_i\} = \{(1, s_i,\ldots,s_i^{d-1}) \}_{i=1}^n$$ 
in ${\mathbb R}^d$, where $s_1 < s_2 < \cdots < s_n$,
i.e. $Z(n,d)$ is the set of all linear combinations $\sum_i c_i u_i$
with $0 \leq c_i \leq 1$.

\smallskip
Let $\Lambda_n = \{0, +, -\}^n$ and $\lambda = 
(\lambda_1, \lambda_2,\ldots,\lambda_n) \in \Lambda_n$. 
An \textit{even gap} of $\lambda$ is a pair of indices 
$i<j$ such that $\lambda_i, \lambda_j$ are nonzero 
entries of \textit{opposite} sign which are separated by an 
even number of zeros, i.e. $\lambda_r = 0$ for all $i < r < 
j$ and $j-i-1$ is even, possibly zero. An \textit{odd gap} 
is a pair of indices $i<j$ such that $\lambda_i, \lambda_j$ 
are nonzero entries of the \textit{same} sign which are 
separated by an odd number of zeros. We define the 
quantity $m(\lambda)$ to be the sum of the number of even 
gaps, the number of odd gaps and the number of zero 
entries of $\lambda$. For example, if $\lambda = 
(+, +, 0, -, -, 0, -, -, +, +, -)$ then $m(\lambda) = 2+1+2 = 
5$, accounting for empty gaps as well.     

Partially order the set $\Lambda_n$ by extending 
componentwise the partial order on $\{0, +, -\}$ defined 
by the relations $+ < 0$ and $- < 0$. 
\begin{proposition} \label{zonclic}
The poset of proper faces of $Z(n,d)$ is isomorphic to the 
induced subposet of $\Lambda_n$ which consists of the 
$n$-tuples $\lambda$ satisfying $m(\lambda) \leq d-1$.
\end{proposition} 
\begin{proof}
Recall that the face poset of $Z(n,d)$
is anti-isomorphic to the poset of 
covectors of the polar hyperplane arrangement 
\cite[\S 7.3]{Ziegler}. These covectors are all possible 
$n$-tuples of the form
\begin{equation}
\lambda = (\sign f(s_1), \sign f(s_2),\ldots,\sign f(s_n)),
\label{sgnf}
\end{equation}
where $f$ is a polynomial of degree at most $d-1$.
It follows from elementary properties of 
polynomials that $f$ has at least $m(\lambda)$ zeros, 
counting multiplicities, so $m(\lambda) \leq d-1$ unless
$\lambda$ is the zero vector. Conversely, given 
$\lambda$, one can construct a polynomial $f$ of degree 
$m(\lambda)$ that satisfies (\ref{sgnf}) by locating its 
zeros and choosing the sign of the leading coefficient 
appropriately.   
\end{proof}

We now turn to the combinatorics of the monotone path 
polytope $\Sigma(C(n,d)\overset{\pi}{\rightarrow} C(n,1))$. Recall
>from Theorem 2.2 that the face poset of 
$\Sigma(C(n,d)\overset{\pi}{\rightarrow} C(n,1))$ is isomorphic to 
the poset of $\pi$-coherent subdivisions of $C(n,1)$. 
The $\pi$-induced subdivisions in this case correspond to the
\textit{cellular strings} \cite{BKS} on $C(n,d)$ with 
respect to $\pi$. These are sequences $\sigma 
= (F_1, F_2,\ldots,F_k)$ of faces of $C(n,d)$ having
the property that $v_1 \in F_1, v_n \in F_k$ and
$\max(\pi(F_i)) = \min(\pi(F_{i+1}))$ for $1 \leq i <k$.
Such a $\sigma$ gives rise to a vector $\lambda = \lambda_{\sigma} = 
(\lambda_2,\ldots,\lambda_{n-1}) \in \Lambda_{n-2}$ as 
follows: For $2 \leq i \leq n-1$ let $\lambda_i = +$, 
respectively $-$, if the vertex $v_i$ of $C(n,d)$ does not 
appear in $\sigma$, respectively is an initial or terminal
vertex, with respect to $\pi$, of some face 
$F_r$ of $\sigma$ and $\lambda_i = 0$ otherwise. For 
example, if $n = 10$ and the faces of $\sigma$ have vertex 
sets $\{v_1, v_3, v_4\}$, $\{v_4, v_7\}$ and $\{v_7, v_8, 
v_{10}\}$ then 
$\lambda _{\sigma} =
 (\lambda_2,\ldots,\lambda_{n-1})= (+, 0, -, +, +, -, 0,+)$.

\smallskip
Recall that $\sigma_1 \leq \sigma_2$ in the Baues poset if and only if
the union of the faces of $\sigma_1$ is contained in the union of the
faces of $\sigma_2$. For cellular strings on $C(n,d)$, this happens 
if and only if $\lambda_{\sigma_1} \leq 
\lambda_{\sigma_2}$ in $\Lambda_{n-2}$. It follows that 
the face poset of $\Sigma(C(n,d)\overset{\pi}{\rightarrow} C(n,1))$ is 
isomorphic to the induced 
subposet of $\Lambda_{n-2}$ which consists of the tuples 
of the form $\lambda_{\sigma}$ for all coherent cellular 
strings $\sigma$ on $C(n,d)$.

\smallskip
The following lemma will be used in the proof of Theorem \ref{monopp},
and is closely related to Lemma 2.3 of \cite{AmeZie}
(incidentally, \cite{AmeZie} also contains very interesting
enumerative aspects of projections of polytopes polar to the cyclic polytopes). 
\begin{lemma} \label{isomorphism}
Let $P, P^{\prime}$ be polytopes with face posets $L,
L^{\prime}$ respectively. Suppose that $\dim(P) \geq 
\dim(P^{\prime})$ and that $\phi: L \toto L^{\prime}$
satisfies $\phi(x) \leq \phi(y)$ if and only if $x \leq y$
for all $x, y \in L$. Then $\phi$ is an isomorphism, i.e. 
a combinatorial equivalence between $P$ and 
$P^{\prime}$. 
\end{lemma}
\begin{proof}
The hypothesis on $\phi$ implies that $\phi$ is injective and sends chains
of $L$ to chains of $L'$. Thus, $\phi$ induces a simplicial injective map
>from the order complex of $L$ into that of $L'$. These order complexes are
isomorphic to the barycentric subdivisions of the boundary complexes of the
polytopes $P$ and $P'$, respectively. Injectivity then implies
that $dim(P)\le dim(P')$. Since a topological sphere cannot
properly contain another topological sphere of the same dimension
(see e.g. \cite[Theorem 6.6. and Exercise 6.9, pp. 67-68]{Massey}),
the simplicial map is bijective, hence an isomorphism of simplicial
complexes, and hence $\phi$ is an isomorphism of posets.
\end{proof}

\noindent
\textit{Proof of Theorem \ref{monopp}}.
Recall that
the face posets of both $\Sigma(C(n,d)\overset{\pi}{\rightarrow} C(n,1))$
and
$Z(n-2,d-1)$
are isomorphic to certain induced subposets of
$\Lambda_{n-2}$. It suffices to show that if $\sigma$
is a coherent cellular string on $C(n,d)$ with respect
to $\pi$ then $m(\lambda_{\sigma}) \leq d-2$.
Indeed, it then follows that there is a well defined map
$\phi = \phi_{n,d}$ from the face poset of
$\Sigma(C(n,d)\overset{\pi}{\rightarrow} C(n,1))$
to that of $Z(n-2,d-1)$ that satisfies the hypothesis of
Lemma \ref{isomorphism}. The face of
$\Sigma(C(n,d)\overset{\pi}{\rightarrow}
C(n,1))$ defined by the
coherent cellular string $\sigma$ is mapped under
$\phi$ to the face of $Z(n-2,d-1)$ which corresponds
to $\lambda_{\sigma}$ under the isomorphism of
Proposition 2.1. Since both polytopes have dimension
$d-1$, the lemma completes the proof.

So suppose that $\sigma$ is coherent and let
$\lambda_{\sigma} = \lambda =
(\lambda_2,\ldots,\lambda_{n-1})$, $\lambda_1 =
\lambda_n = -$. By Ziegler's definition of
$\pi$-coherence \cite[\S 9.1]{Ziegler},
there is a polynomial $f$ of degree at most $d$ such that the
polygon $\hat{Q}_f := \conv \{(t_i, f(t_i)): \ 1 \leq i \leq n\}$
has the following property: the points $(t_i, f(t_i))$ for
which $\lambda_i = -$ are the lower vertices of $\hat{Q}_f$,
the ones for which $\lambda_i = 0$ lie on its lower
edges and the ones for which $\lambda_i = +$ lie
above. In the rest of the proof we show that $f$ has degree at
least $m(\lambda)+2$, so that $m(\lambda) \le d-2$.

We assume that there is no $2 \leq i < n-1$ such that $\lambda_{i} =
\lambda_{i+1} = +$, since otherwise we can drop any of the two indices 
$i$ or ${i+1}$ without decreasing the value of $m(\lambda)$.
Let $l_i$ be the line segment joining the points
$(t_{i},f(t_{i}))$ and $(t_{i+1},f(t_{i+1}))$, for $1 \le i < n$. We
construct
the $(n-2)$-tuple $\lambda''=(\lambda''_2,\dots,\lambda''_{n-1})$ as
follows: $\lambda''_i$ equals $+$, $0$ or $-$ depending on whether the
slope of the segment $l_i$ is smaller, equal or greater than the slope
of the segment $l_{i-1}$.

It is easy to verify that $\lambda_i\in\{+,-\}$
implies $\lambda''_i=\lambda_i$. On the other hand, $\lambda_i=0$
implies that $\lambda''_i=0$ unless at least one of $\lambda_{i-1}$ or
$\lambda_{i+1}$ equals $+$, in which case $\lambda''_i=-$. In other
words, the $(n-2)$-tuple $\lambda''$ is obtained from $\lambda$ by
changing every pair of consecutive entries $(+, 0)$ to $(+, -)$ and
every pair $(0, +)$ to $(-, +)$. This implies that $m(\lambda'')
= m(\lambda)$. For $1 \leq i < n$ let $\theta_i$
be such that $t_{i} < \theta_i < t_{i+1}$ and
\[ \frac{f(t_{i+1}) - f(t_{i})}{t_{i+1} - t_{i}} =
f^{\prime}(\theta_i). \]
Observe that $f^{\prime}(\theta_i)$ equals the slope of the segment
$l_i$. For $2\leq i \leq n-1$ let $\mu_i$
be such that $\theta_{i-1} < \mu_i < \theta_{i}$ and
\[ \frac{f'(\theta_{i}) - f'(\theta_{i-1})}{\theta_{i} - \theta_{i-1}} =
f''(\mu_i). \]
Then $ \sign f''(\mu_i) = \sign (f'(\theta_{i}) -
f'(\theta_{i-1})) = - \lambda''_i$. This implies that the degree of
$f''$ is at least $m(\lambda'')$ and finishes the proof.
\qed

The combinatorial equivalence $\phi_{n,d}$, described in the
proof of Theorem \ref{monopp}, together with Proposition \ref{zonclic} 
implies the following corollary.
\begin{corollary}
\label{coro:mlambda}A cellular string $\sigma$ for 
$\pi:C(n,d) \rightarrow C(n,1)$ is $\pi$-coherent 
if and only if $m(\lambda_{\sigma}) \leq d-2$. In particular, 
whether $\sigma$ is coherent or not does not depend on the 
choice of $t_1,\ldots,t_n$ used to define $C(n,d)$.
\end{corollary}
A cellular string $\sigma = (F_1, F_2,\ldots,F_k)$ on $C(n,d)$ 
with respect to $\pi$ is a tight $\pi$-induced 
subdivision if all its faces $F_r$ are edges $v_{i_{r-1}} 
v_{i_r}$ of $C(n,d)$, where $1 = i_0 < i_1 < \cdots < i_k = n$. 
Equivalently, $\sigma$ is tight if $\lambda_{\sigma}$ 
contains no zeros. The vertices $v_i$ of $\sigma$ correspond 
to the indices $i$ for which $\lambda_i = -$,
together with the indices $1$ and $n$. The tight cellular 
strings are the \textit{monotone edge paths}. 
 
 As a corollary, we characterize and enumerate the monotone 
edge paths on $C(n,d)$ which are $\pi$-coherent
for $\pi:C(n,d)\rightarrow C(n,1)$.
For $\lambda \in \{+, -\}^{n-2}$ let $c(\lambda)$ be the number of 
maximal strings of successive $+$ signs or successive 
$-$ signs in $\lambda$. Note that $m(\lambda) = c(\lambda) 
- 1$. 
\begin{corollary}
For $\lambda \in \{+, -\}^{n-2}$ we have $\lambda = 
\lambda_{\sigma}$ for a $\pi$-coherent monotone 
edge path $\sigma$ in $C(n,d)$ if and only if $c(\lambda) \leq 
d-1$. The number of $\pi$-coherent monotone edge 
paths in $C(n,d)$ is 
\[ 2 \sum_{j=0}^{d-2} {n-3 \choose j}.\]
\end{corollary}
\begin{proof}
The first statement follows from Corollary \ref{coro:mlambda}. For the 
second statement, note that there are $2 {n-3 \choose j}$
tuples $\lambda \in \{+, -\}^{n-2}$ with $c(\lambda) 
= j+1$. 
\end{proof}
\begin{remark} \rm \ \\ 
\noindent
If $d \geq 4$, any two vertices of $C(n,d)$ are connected 
by an edge and the total number of monotone edge paths 
on $C(n,d)$ is $2^{n-2}$. Hence when $d$ is fixed and $n$ gets
large, the fraction of coherent paths approaches zero.
Similar behavior is exhibited in
Proposition 5.10 of \cite{dhss}, where it is proved that the cyclic
polytope $C(n,n-4)$ has $\Omega(2^n)$ triangulations but only $O(n^4)$
regular ones.  Another example of this behavior with regard to monotone
paths for non-cyclic polytopes appears in \cite{Ath}. 
\end{remark}

\begin{remark} \rm \ \\
\noindent
One can rephrase Theorem \ref{monopp} as saying that for the linear
functional $f({\bf x})=x_1$ mapping $C(n,d)$ onto a $1$-dimensional
polytope $f(C(n,d))$, the monotone path polytope
$$\Sigma(C(n,d) \overset{f}{\rightarrow} f(C(n,d)))$$ has face lattice
independent of the choice of $t_1,\ldots,t_n$. One might ask whether
this is true for all linear functionals. It turns out that this is not
the case.
Consider the linear functional $f({\bf x})=x_1+x_3$ on the polytope
$C(5,3)$. Since the functional is monotone along the moment curve, it will 
produce the same monotone paths (actually the same ones as the
standard functional $x_1$) for any choice of parameters $t_1<\cdots<t_5$. 
However, different choices of parameters can change the set of coherent 
monotone edge paths.
\end{remark}

Let us fix $t_3=0$ and $t_2=-t_4$, $t_1=-t_5$, so that we only have
two free parameters $0<t_4<t_5$ and the cyclic polytope $C(5,3)$
has a symmetry $[x_1\to -x_1; x_3\to -x_3]$ which exchanges the
vertices $i$ and $6-i$. We want to find out when the monotone
path consisting of the edges $13$ and $35$ is coherent. The normal
vectors to the faces $235$ and $345$ are $(t_4t_5,-t_4-t_5,1)$ and 
$(-t_4t_5,t_4-t_5,1)$, respectively, so these two vectors
generate the two boundary rays in the normal cone to the segment
$35$. Thus, the projection of the normal cone of $35$ to 
$ker(f)$ is bounded by the vectors
$({t_4t_5-1\over 2},-t_4-t_5,{1-t_4t_5\over 2})$ 
and $({-t_4t_5-1\over 2},t_4-t_5,{1+t_4t_5\over 2})$. By symmetry, the
projection to $ker(f)$ of the normal cone to the segment $13$ is
bounded by $({1-t_4t_5\over 2},-t_4-t_5,{t_4t_5-1\over 2})$
and $({t_4t_5+1\over 2},t_4-t_5,{-1-t_4t_5\over 2})$. If $t_4t_5>1$
then the relative interiors of these two cones in $ker(f)$ intersect
(in the vector $(0,-1,0)$ for example). If $t_4t_5\le 1$ then they do
not intersect since the first one only contains vectors with $x_1<0$
and the second one vectors with $x_1>0$.
Thus, the path containing the segments $13$ and $35$ is coherent if
and only if $t_4t_5>1$.

\begin{remark} \rm \ \\
\noindent
Recall that the Upper Bound Theorem \cite[\S 8.4]{Ziegler},
\cite[\S 4.7]{Grunbaum} states that the cyclic polytope $C(n,d)$ has
the most boundary $i$-faces among all $d$-polytopes with $n$ vertices
for all $i$.  We have also seen that the facial structure
of the monotone path polytope 
$\Sigma(C(n,d) \overset{\pi}{\rightarrow} C(n,1))$
is independent of the choice of points on the moment curve.
These two facts might tempt one to make the following
``Upper Bound Conjecture (UBC)
for monotone path polytopes": For all $d$-polytopes with $n$-vertices
and linear functionals $f$, the monotone path polytope 
$\Sigma(P \overset{f}{\rightarrow} f(P))$
has no more boundary $i$-faces than 
$\Sigma(C(n,d) \overset{\pi}{\rightarrow} C(n,1))$.
However, this turns out to be false, as demonstrated by the
example of a non-neighborly simplicial 4-polytope with eight vertices whose 
monotone path
polytope with respect to the projection to the first coordinate
has two more coherent paths than C(8,4). The vertex coordinates of this
$4$-polytope are given by the columns of the following matrix:

$$\left [\begin {array}{cccccccc} -84&-36&-35&11&90&31&47&-50
\\\noalign{\medskip}-54&71&-71&-17&65&-34&60&99\\\noalign{\medskip}48&
36&73&-40&50&54&24&65\\\noalign{\medskip}6&-65&52&100&-39&49&-76&-15
\end {array}\right ]$$

\noindent
This raises the following question.
\begin{question}
Is there some natural family of polytope projections $P
\overset{\pi}{\rightarrow} Q$ indexed by $(n,d',d)$ with $dim(P)=d',
dim(Q)=d$, such that $P$ has $n$ vertices and the
fiber polytope $\Sigma(P \overset{\pi}{\rightarrow} Q)$ has more
$i$-faces than any other fiber polytope of an $n$-vertex $d'$-polytope
projecting onto a $d$-polytope?
\end{question}
For the case $d=0$, the Upper Bound Theorem says that
the family of projections $C(n,d') \overset{\pi}{\rightarrow} C(n,d)$
provides the answer, but the above counterexample shows that it
does not already for $d=1$. However, one could still
ask whether this family provides the answer asymptotically. For simplicity,
we restrict our attention to the case $d=1$ of monotone path polytopes
and the number of vertices of $\Sigma(P \overset{\pi}{\rightarrow} Q)$. 
We also replace $d'$ by $d$, as in the beginning of this section. If $d$
is fixed, Corollary 3.5 implies that the number of vertices of 
$\Sigma(C(n,d) \overset{\pi}{\rightarrow} C(n,1))$ is a polynomial in
$n$ of degree $d-2$. 

Let $r_d (n)$ denote the maximum
number of vertices that a monotone path polytope of an $n$-vertex 
$d$-polytope projecting onto a line segments can have.
\begin{question}
With $d$ fixed, what is the 
asymptotic behaviour of $r_d (n)$ as $n \rightarrow \infty$? In particular, 
is $r_d (n)$ bounded above by a polynomial in $n$ of degree $d-2$?
\end{question}
\end{remark}
We close this section by giving a polynomial upper bound for $r_d (n)$
of degree $3d-6$. Let 
\[ q_d (m) = \sum_{j=0}^{d} {m \choose j}. \]
\begin{proposition}
We have $r_d (n) \leq 2 \, q_{d-2} ({n \choose 3}-1)$, a polynomial in $n$ of 
degree $3d-6$.
\end{proposition}
\begin{proof}
Let $P \subseteq {\mathbb R}^d$ be $d$-dimensional with vertices $p_1, 
p_2,\dots,p_n$ and $\pi : {\mathbb R}^d \rightarrow {\mathbb R}$ be a linear 
map. Let $a_i = \pi(p_i)$ for $1 \leq i \leq n$. 

We use Ziegler's definition \cite[\S 9.1]{Ziegler} of a $\pi$-coherent
monotone edge path, described in Section \ref{background}. Let $f \in
ker(\pi)^*$ be a generic linear functional and $\hat{Q}$ be the
convex hull of the points $\hat{a}_i := (a_i, f(p_i))$ in ${\mathbb R}^2$, 
for $1 \leq i \leq n$. The set of lower edges of $\hat{Q}$ is determined by
the oriented matroid of the point configuration ${\mathcal A} :=
\{\hat{a}_i\}_{i=1}^n$. Equivalently, it is determined by the data which
record for each triple $1 \leq i < j < k \leq n$ which of the two 
halfplanes determined by the line though $\hat{a}_i$ and $\hat{a}_k$ the 
point $\hat{a}_j$ lies on. This is equivalent to recording which side of a 
certain linear hyperplane in $ker(\pi)^*$, depending on $(i,j,k)$, 
the functional $f$ lies on. Hence the number of $\pi$-coherent
monotone edge paths on $P$ is at most the number of regions
into which some ${n \choose 3}$ linear hyperplanes dissect 
$ker(\pi)^*$. This number is at most the number of regions into which
a \textit{generic} arrangement of ${n \choose 3}$ linear hyperplanes
dissects ${\mathbb R}^{d-1}$, which is $2 \, q_{d-2} ({n \choose 3}-1)$
(see \cite{Zaslavsky} for more on
counting regions in hyperplane arrangements).
\end{proof}

\section{The case $d'=n-1$: triangulations and secondary polytopes}
\label{secondary-polytopes}

In this section we restrict our attention to the natural  
projection $\pi: C(n,n-1) \to C(n,d)$ and prove the assertions
of Theorem \ref{main-theorem} concerning the case $d'=n-1$.
In this case, since $C(n,n-1)$ is an $(n-1)$-simplex $\Delta^{n-1}$
(and since all $(n-1)$-simplices are affinely equivalent), the fiber polytope
$\Sigma(C(n,n-1)\overset{\pi}{\rightarrow}C(n,d))$ coincides with
the \textit{secondary polytope} $\Sigma(C(n,d))$
(see Section \ref{background} or \cite[\S 7]{GKZbook}),
whose vertices correspond to the regular triangulations of $C(n,d)$.
The question of the existence of non-regular triangulations of $C(n,d)$
was first raised in \cite[Remark 3.5]{KV}. Billera, Gel'fand and
Sturmfels first constructed such a triangulation for $C(12,8)$ in
\cite[\S 4]{BGS}.  Our results show that this example is far from minimal
and provide a complete characterization of the values of $n$ and $d$ for
which $C(n,d)$ has non-regular triangulations.

We recall and separate out the assertions of Theorem \ref{main-theorem}
which deal with secondary polytopes.  In this context, we use the terms
``coherent subdivision" and ``regular subdivision" interchangeably, as
both occur in the literature. 

\begin{theorem}
\label{seccyc}
All polytopal subdivisions of $C(n,d)$ are coherent if and only if either
\begin{enumerate}
\item[$\bullet$] $d \leq 2$ or
\item[$\bullet$] $n-d \leq 3$ or
\item[$\bullet$] $(n,d) \in \{(8,4),(8,3),(7,3)\}$.
\end{enumerate}
In all other cases, there exists a subdivision of $C(n,d)$ whose
coherence varies with the choice of points on the moment curve,
and for every choice of points there is some
incoherent subdivision.
\end{theorem}

The proof of this result occupies the remainder of this section.
We begin by showing that in all of the cases asserted above, all
polytopal subdivisions are coherent.  For $d = 1$ this is easy
and for $d=2$ and $n-d \leq 3$ this was shown by Lee \cite{Lee1, Lee2}.
In fact, these references show that all subdivisions in these cases
are {\it placing} subdivisions (see definition in \cite{Lee2}).
It therefore remains to show that all subdivisions of $C(n,d)$ are coherent
for $(n,d)$ equal to $(8,4),(8,3)$ or $(7,3)$.

Our task is simplified somewhat by the following fact.

\begin{lemma}  \label{oneall}
Suppose that for a certain choice of points along the moment curve,
the canonical projection $C(n,d') \overset{\pi}{\rightarrow} C(n,d)$
has the property that every $\pi$-induced
subdivision is $\pi$-coherent. Then, the same happens for every other choice
of points along the moment curve.
\end{lemma}

\begin{proof}
In every choice of points the poset of $\pi$-coherent subdivisions is
isomorphic to the face poset of the corresponding fiber polytope,
which is a polytope of dimension $d'-d$. 
The hypothesis of the lemma implies that, for a certain
choice of points, the face poset of the fiber polytope is the whole Baues
poset $\Baues(C(n,d') \overset{\pi}{\rightarrow} C(n,d))$.
Since the Baues poset is independent of the choice of points,
Lemma 4.4 of \cite{BS2} implies 
that in every other choice of points the Baues poset coincides with the
face poset of the fiber polytope.
\end{proof}

\begin{remark} \rm \ \\
\noindent Lemma \ref{oneall} is true in a more general
situation, namely, whenever we have two projections of polytopes
$P \overset{\pi}{\rightarrow} Q$ and $P' \overset{\pi'}{\rightarrow}
Q'$ and there is a bijection $\phi$ between the vertices of $P$ and $P'$ which
induces an isomorphism between the oriented matroids of affine
dependencies of $P$ and $P'$, as well as those of $Q$ and $Q'$.
This is so because these assumptions imply that the two Baues posets are
isomorphic.

On the other hand, it is not enough to assume that $\phi$ induces only a
combinatorial equivalence for $P$ and $P'$ and for $Q$ and $Q'$.
For example, if $P$ and $P'$ are two
5-simplices projecting in the natural way onto two combinatorial
octahedra $Q$ and $Q'$ with different oriented matroid (i.e.
different affine dependence structure) then both Baues posets contain
all the polytopal subdivisions of $Q$ and $Q'$ respectively, but they
are different and the proof of Lemma \ref{oneall} is not valid.
This is relevant to the situation with cyclic polytopes $C(n,d)$ since
there exist polytopes with the same face lattice as $C(n,d)$ but whose
vertices have different affine dependence structure \cite{Bisztriczky-Karolyi}.
\end{remark}

\begin{remark} \rm \ \\
We also note that Lemma \ref{oneall} shows that the last assertions
in parts (2) and (4) of Theorem \ref{main-theorem} follow from the
assertions preceding them.  To be precise, if $(n,d,d')$ are such
that there exists some $\pi$-induced subdivision of $C(n,d)$
whose $\pi$-coherence depends upon the choice of points
on the moment curve, then Lemma \ref{oneall} implies
that there cannot exist a choice of points for which every
$\pi$-induced subdivision of $C(n,d)$ is $\pi$-coherent.  
\end{remark}

\vskip 0.1 in

It was recently shown by Rambau \cite{Rambau1} that all triangulations of 
$C(n,d)$ are connected by bistellar flips. Hence one 
can rely on this fact to enumerate all triangulations in small instances (see
Table \ref{tricnd}). The program PUNTOS is
an implementation of this algorithmic procedure and can be obtained via 
anonymous ftp at {\tt ftp://geom.umn.edu, directory /priv/deloera} (see 
\cite{thesjes} for details). In Lemma \ref{subdivisions}
we will use the information given by PUNTOS
for the three cases which interest us to prove 
that all the subdivisions are regular in a certain
choice of points along the moment curve, and hence in all choices by
Lemma \ref{oneall}. 

The following lemma is a direct proof of the fact that all the
triangulations are regular, in every choice of points, for the three
cases. The lemma also clearly follows from Lemma \ref{subdivisions}, but we
find the proof below of independent interest.

\begin{lemma} All triangulations of $C(7,3)$ and $C(8,4)$ are placing.
All triangulations of $C(8,3)$ are regular.
 \label{cndplace} 
\end{lemma}

\begin{proof}
We know from the results of \cite{Lee2} that all triangulations of $C(6,3)$ and
$C(7,4)$ are placing. It is enough to check that each triangulation
of $C(7,3)$ and $C(8,4)$ has at least one vertex whose link is contained in
the boundary complex of $C(6,3)$ and $C(7,4)$ respectively. 
The polytope $C(7,3)$ has a single symmetry that maps $i$ to $7-i+1$.
The original 25 triangulations are divided into 16 distinct symmetry 
classes. For $C(8,4)$, which has the dihedral group of order 16 as
its group of symmetries, the original 40 triangulations
are divided into only 4 symmetry classes.
In Tables \ref{markac73} and \ref{markac84} we show, for each of the 16 
symmetrically distinct triangulations of $C(7,3)$ and each of the 4
 symmetrically distinct triangulations of $C(8,4)$, that the above 
``good link'' property is indeed satisfied at least at one vertex. 

{\footnotesize
\begin{table}[tb]
\begin{center}
\begin{tabular}{| l| l |} \hline
 triangulations of $C(7,3)$ modulo symmetries & good \\ 
                           & links at:\\ \hline \hline

2356,1234,4567,3467,2345,2367,1256,3456,1267,1245& 1,7\\ \hline 
2456,2346,1234,4567,3467,2367,1256,1267,1245&1,7\\ \hline 
2356,1234,2345,2367,1256,1267,1245,3457,3567&1\\ \hline 
2346,1234,4567,3467,2367,1267,1456,1246&5,7\\ \hline 
2356,2367,1256,1267,1235,1345,3457,3567&4\\ \hline 
1234,2345,1256,1267,1245,3457,2567,2357&1\\ \hline 
2367,1267,1345,3457,3567,1236,1356&4\\ \hline 
4567,3467,3456,1345,1356,1237,1367&2\\ \hline 
4567,3467,2367,1267,1456,1236,1346&5,6,7\\ \hline 
4567,3467,1456,1237,1367,1346&2,5\\ \hline 
1345,3457,3567,1356,1237,1367&2,3,4\\ \hline 
1345,3457,1237,1357,1567&2,4,6\\ \hline 
1237,1567,1457,1347&1,2,6,7   \\ \hline 
1234,2347,1567,1457,1247&3,6\\ \hline 
1234,4567,1456,2347,1247,1467&3,4,5\\ \hline 
1234,4567,1267,1456,1246,2347,2467&3,5\\ \hline 
\end{tabular}
\caption{Triangulations of $C(7,3)$ are placing triangulations}
 \label{markac73}
\end{center}
\end{table}
}
{\footnotesize
\begin{table}[ptb]
\begin{center}
\begin{tabular}{| l| l |} \hline
 triangulations of $C(8,4)$ modulo symmetries  & good \\ 
                           & links at:\\ \hline \hline

23678,23458,12568,12458,45678,23568,12678,12348,34568,34678 & 1,7,8\\ \hline
24568,23456,23678,12568,12458,45678,12678,12348,34678,23468 & 1,7 \\ \hline
23678,12568,45678,23568,12678,34568,34678,13458,12358,12345 & 7 \\ \hline
23678,45678,12678,34568,34678,13458,12345,12368,12356,13568  &  7\\ \hline
\end{tabular}
\caption{Triangulations of $C(8,4)$ are placing triangulations}
 \label{markac84}
\end{center}
\end{table}
}

In the case of $C(8,3)$ one can verify that, modulo symmetries, only
the following five triangulations are not placing triangulations
(if symmetries are not considered, these are 8 triangulations out of
a total of 138):

\vskip 12pt
\begin{obeylines}
$2378,2356,2367,1267,3456,3478,3467,1256,1278,1345,1235,4568,4678$
$2378,2367,1267,3456,3478,3467,1278,1345,4568,4678,1236,1356$
$2356,1267,3456,1256,1278,1345,1235,4568,3468,2678,2368$
$1267,3456,1278,1345,4568,1236,1356,3468,2678,2368$
$2378,2367,1267,3456,1278,1345,4568,1236,1356,3678,3468$.
\end{obeylines}
\vskip 12pt

One can check that each of the first four triangulations has a neighbor
which is a placing triangulation having only 4 bistellar flips. This
implies that if any of the four were to disappear from the
4-dimensional secondary polytope $\Sigma(C(8,3))$ for some choice of
points on the moment curve then the placing triangulation in question
would be left with three neighbors, which is impossible in a
4-dimensional polytope.

This leaves only the fifth triangulation on the above list whose
regularity must be checked directly. 
Regularity can be determined via the feasibility of a certain system of
linear inequalities. The variables of this system are the heights $w_i$. 
The inequalities are determined by pairs of maximum dimensional simplices 
and points. Each inequality establishes the fact that a given point
lies ``above'' a certain hyperplane after using the heights $w_i$ for
a lifting. The coefficients of the inequalities can be interpreted as 
oriented volumes. These inequalities form a system $B \bar{\lambda} <0$ 
which will be feasible precisely when the triangulation $K$ is regular.

Farkas' theorem \cite{schrijver} indicates that exactly one of the 
following holds: either the system of inequalities $B \bar{\lambda} 
<0 $ defined by the triangulation is consistent, or there exists $y
\in {\mathbb R}^m$ such that $y^{T}B=0, \quad \sum y_i>0$ and $y \geq 0$. 
The support of such a vector $y$ labels an inconsistent subset of
inequalities in $B \bar{\lambda} <0$. An explicit proof of non-regularity 
is then the impossible inequality $0=(y^{T}B)\bar{\lambda} =y^{T} (B 
\bar{\lambda}) <0$. Notice that such a vector $y$ lies in the
kernel of the transpose of $B$.

One can set up the following matrix $B$ for the
last triangulation of the list. The coefficients are simply 
Vandermonde determinants because they determine the volume of the
simplices in the triangulation.

{\tiny
$$\left [\begin {array}{cccccccc} -{\it vol}(2367)&{\it vol}(1367)&-{
\it vol}(1267)&0&0&{\it vol}(1237)&-{\it vol}(1236)&0
\\\noalign{\medskip}0&0&-{\it vol}(4568)&{\it vol}(3568)&-{\it vol}(
3468)&{\it vol}(3458)&0&-{\it vol}(3456)\\\noalign{\medskip}0&-{\it 
vol}(3678)&{\it vol}(2678)&0&0&-{\it vol}(2378)&{\it vol}(2368)&-{\it 
vol}(2367)\\\noalign{\medskip}-{\it vol}(2378)&{\it vol}(1378)&-{\it 
vol}(1278)&0&0&0&{\it vol}(1238)&-{\it vol}(1237)\\\noalign{\medskip}-
{\it vol}(2678)&{\it vol}(1678)&0&0&0&-{\it vol}(1278)&{\it vol}(1268)
&-{\it vol}(1267)\\\noalign{\medskip}{\it vol}(2356)&-{\it vol}(1356)&
{\it vol}(1256)&0&-{\it vol}(1236)&{\it vol}(1235)&0&0
\\\noalign{\medskip}-{\it vol}(3456)&0&{\it vol}(1456)&-{\it vol}(1356
)&{\it vol}(1346)&-{\it vol}(1345)&0&0\\\noalign{\medskip}0&0&{\it vol
}(4678)&-{\it vol}(3678)&0&{\it vol}(3478)&-{\it vol}(3468)&{\it vol}(
3467)\end {array}\right ]$$}

For example, the first row of the matrix $B$ corresponds to the
inequality that indicates that point 1 is above the lifted plane $2367$.
The kernel of the transpose of the above matrix is 4-dimensional and
spanned by the columns of the following matrix, where $t_{ij}$ denotes
the difference $t_i-t_j$ for $i<j$.

$$
\left [
\begin {array}{cccc}
-{\frac {(t_{78}t_{38}t_{28})}{(t_{67}t_{36}t_{26})}}&
-{\frac {(t_{78}t_{68}t_{28})}{(t_{37}t_{36}t_{23})}}&
{\frac {(t_{56}t_{35}t_{25})}{(t_{67}t_{37}t_{27})}}&
-{\frac {(t_{56}t_{35}t_{46}t_{45} t_{34})}
              {(t_{67}t_{37}t_{27}t_{26}t_{23})}}\\
\noalign{\medskip}
0&0&
-{\frac {(t_{26}t_{23}t_{16}t_{13}t_{12})}
              {(t_{68}t_{38}t_{48}t_{46}t_{34})}}&
{\frac {( t_{16}t_{13}t_{14})}{(t_{68}t_{38}t_{48})}}\\
\noalign{\medskip}
-{\frac {(t_{17}t_{13}t_{12})}{(t_{67}t_{36}t_{26})}}&
-{\frac {(t_{17}t_{16}t_{12})}{(t_{37}t_{36}t_{23})}}&
-{\frac {(t_{56}t_{35}t_{16}t_{13}t_{57}t_{12})}
              {(t_{78}t_{68}t_{67}t_{38}t_{37}t_{27})}}&
-{\frac {( t_{56}t_{35}t_{46}t_{45}t_{34}t_{16}t_{13}t_{17})}
              {(t_{78}t_{68}t_{67}t_{38}t_{37}t_{27}t_{26}t_{23})}}\\
\noalign{\medskip}
1&0&0&0 \\
\noalign{\medskip}
0&1&0&0 \\
\noalign{\medskip}
0&0&1&0 \\
\noalign{\medskip}
0&0&0&1 \\
\noalign{\medskip}
0&0&
-{\frac {( t_{56}t_{35}t_{26}t_{23}t_{16}t_{13}t_{12}t_{58})}
              {(t_{78}t_{68}t_{67}t_{38}t_{37}t_{48}t_{46}t_{34})}} &
-{\frac {( t_{56}t_{35}t_{16}t_{13}t_{18}t_{45})}
              {(t_{78}t_{68}t_{67}t_{38}t_{37}t_{48})}}\\
\end {array}
\right ]
$$

The third row of the above matrix is negative regardless of the values of
the parameters $t_i$. Since this matrix also contains a $4 \times
4$ identity submatrix, there is no non-zero vector $y$ in its column
space that has all its elements nonnegative. The system defined by the
matrix $B$ is always feasible, making the fifth triangulation regular
regardless of the values of $t_i$.
\end{proof}

We define the \textit{ranking} of a polytopal subdivision of $C(n,d)$
to be the sum of the dimensions of the secondary polytopes of their
disjoint cells that are
not simplices. For instance, for $C(8,4)$ it is possible to have a polytopal
subdivision with two copies of $C(6,4)$ as the only cells that are not 
simplices; the ranking of such subdivisions is 2.  We do
not use the word rank as the poset of polytopal subdivisions of $C(n,d)$ is
not necessarily graded. 

\begin{lemma}
\label{subdivisions} 
All polytopal subdivisions of $C(7,3)$, $C(8,3)$ and $C(8,4)$ are
coherent for every choice of parameters.
 \label{allcoherent} 
\end{lemma}

\begin{proof}
Because of Lemma \ref{oneall} we only need to prove that all the
subdivisions are coherent in one choice of parameters. So we fix the
parameters to be $t_i=i$. Also, it suffices to prove the result for
$C(8,3)$ and $C(8,4)$ because $C(7,3)$ is a subpolytope of $C(8,3)$.

We describe a procedure to
enumerate all polytopal subdivisions of $C(n,d)$ of ranking $k$ for
small values of $n$. We say that a polytopal subdivision $S$ of is of 
\textit{type}
$[r_1C(s_1,d),r_2C(s_2,d),\dots,r_mC(s_r,d)]$ if the total number of
cells used in $S$ which are not simplices is $r_1+r_2+\dots+r_m$ and $S$
contains precisely $r_i$ disjoint isomorphic copies of $C(s_i,d)$ with
$s_i>d+1$ ($C(s_i,d)$ is a cell which is not a simplex). For example,
there are 18 polytopal subdivisions of $C(8,4)$ of type
$[2C(6,4)]$. Clearly, all subdivisions of the same type have the same
ranking. The possible isomorphic classes of subpolytopes of
$C(n,d)$ which are not simplices are
$C(d+2,d),C(d+3,d),\dots,C(n-1,d)$. Their secondary polytopes have
dimensions $1,2,\dots ,n-d-2$ respectively. 

Given that we have a complete list of triangulations of $C(n,d)$, to 
count the polytopal subdivisions of $C(n,d)$ of type 
$[r_1C(s_1,d),$ $r_2C(s_2,d), \dots,r_mC(s_m,d)]$ we fix a
triangulation for each $C(s_i,d)$ and form all
possible $(r_1+r_2+\dots+r_m)$-tuples of disjoint triangulated copies of
$C(s_1,d)$, $C(s_2,d),\dots,C(s_m,d)$. Because the copies have been
triangulated, it suffices to count the triangulations of $C(n,d)$
which complete the different tuples of triangulated subpolytopes.
It is easy to list all types of polytopal subdivisions of $C(8,4)$ and
$C(8,3)$ and their cardinalities. 

{\footnotesize
\begin{table}[pb]
\begin{center}
\begin{tabular}{| l| l |} \hline
 type (case of C(8,4)) & cardinality \\ \hline \hline

$[C(7,4)]$ & 8 \\ \hline
$[2C(6,4)]$ & 18 \\ \hline
$[3C(6,4)]$ & 0 \\ \hline \hline

 type (case of C(8,3)) & cardinality \\ \hline \hline

$[2C(5,3)]$  & 162 \\ \hline
$[C(6,3)]$  & 52 \\\hline
$[3C(5,3)]$ & 18 \\\hline
$[C(6,3),C(5,3)]$ & 24 \\  \hline
$[C(7,3)]$ & 8 \\ \hline
$[2C(6,3)]$ & 0 \\ \hline
$[C(6,3),2C(5,3)]$ & 0 \\ \hline
$[4C(5,3)]$ & 0 \\ \hline

\end{tabular}
\caption{Polytopal subdivisions of $C(8,4)$ and $C(8,3)$ by type}
 \label{types8}
\end{center}
\end{table}
}

We computed these numbers using a MAPLE implementation of the above
criteria.  We present the main information in Table \ref{types8}. We
disregard subdivisions of ranking one since these are
exactly the bistellar flips that we can compute
with PUNTOS \cite{thesjes}. Once we get zero subdivisions for all types
of a certain ranking $i$, we do not need to compute the number for any
ranking $j>i$ since any subdivision of ranking $j$ could be refined
into one of ranking $i$.  Also, if a subdivision contains a cell
$C(n-1,d)$ then all the other cells are simplicial and the subdivision
is the one obtained from the trivial subdivision by {\em pushing} the
vertex which is not in the cell $C(n-1,d)$ (see \cite[page
411]{OMbook}). Thus, the number of subdivisions of type $[C(n-1,d)]$
equals $n$ and we do not need to compute any other type containing a
cell $C(n-1,d)$.

Using the program PUNTOS \cite{thesjes} we have computed all
triangulations of $C(8,3)$ and $C(8,4)$ (in the standard choice of
parameters) and checked that they are all regular. This also implies
that all the ranking one subdivisions are regular, since ranking one
subdivisions correspond to bistellar flips between triangulations and
bistellar flips between regular triangulations correspond to edges of
the secondary polytope.
Thus for $C(8,4)$, whose secondary polytope is 3-dimensional, it only
remains to check that the number of subdivisions of ranking at least two
coincides with the number of facets of the secondary polytope. Since
we have 40 triangulations and 64 bistellar flips computed by PUNTOS,
the number of facets of the secondary polytope is given by Euler's
formula and turns out to be 26. This coincides with the results of
Table \ref{types8} and hence all subdivisions are regular.

In the case of $C(8,3)$, PUNTOS tells us that there are 138
triangulations  and 302 ranking one subdivisions
(flips), all of them regular (in the standard choice of parameters).
The calculation of the
secondary polytope using PORTA indicates that for the usual
parameters, the secondary polytope for $C(8,3)$ is a 4-dimensional
polytope which has indeed 50 facets. Then Euler's formula gives a
number of 214 for the number of 2-dimensional faces. Since these
numbers coincide with the numbers of ranking three and ranking two 
subdivisions in Table \ref{types8}, all subdivisions are regular.

\end{proof}

\begin{remark} \rm \ \\
\noindent
In \cite{santos} it is proved that any non-regular subdivision of a polytope
can be refined to a non-regular triangulation. With this, our last
lemma also follows from the fact that in the standard choice of parameters
all triangulations of $C(8,3)$, $C(8,4)$ and $C(7,3)$ are regular,
with no need to analyze non-simplicial subdivisions.
\end{remark}

\medskip

The next lemma shows that in order
to complete the proof of Theorem \ref{seccyc} we need only to exhibit
certain minimal counterexamples.

\begin{lemma} \label{extend} Suppose there exists a 
triangulation $T$ of the cyclic polytope $C(n,d)$ which is regular or
non-regular, depending on the choice of points along the moment
curve. Then such a triangulation also exists for $C(n+1,d)$ and
$C(n+1,d+1)$.
\end{lemma}

\begin{proof}
For the first statement, given such a triangulation $T$ of $C(n,d)$,
extend the triangulation by placing (see definition in \cite{Lee2}) the extra
point $n+1$, that is to say, by joining $n+1$ to all the facets of $T$
which are visible from it. This produces a triangulation $T'$ of $C(n+1,d)$
and it is easy to see that $T'$ is regular for a choice
$t_1 < \cdots < t_n <t_{n+1}$ of parameters on the moment curve
if and only if $T$ is regular for the choice $t_1 < \cdots < t_n$.

For the second statement we use the fact that $C(n+1,d+1)$ is a single
element lifting of $C(n,d)$, so that Lemma \ref{lifting-lemma} applies.
As was discussed in Section \ref{background},
we cannot guarantee that $C(n+1,d+1)$ is a single
element lifting of $C(n,d)$ unless the parameters $t_1 < \cdots
< t_n <t_{n+1}$ are chosen so that $t_{n+1}=0$. However, this presents
no problem since an affine transformation $t_i \mapsto a t_i+b$
of the parameters produces an affine transformation of
the points $v_i$ and thus preserves regularity of triangulations.

By our hypotheses, there exist two choices of parameters for the
cyclic polytope $C(n,d)$ which produce different chamber complexes in
the dual.  By Lemma \ref{lifting-lemma}, the Gale transform of $C(n+1,d+1)$
is obtained from that of $C(n,d)$ by adding a single point.
It is impossible to add a new point to the Gale
transforms and make them equal (in a labeled sense). Thus, there are
choices of parameters for $C(n+1,d+1)$ which produce different chamber
complexes in the Gale transform and which, in particular, produce
different collections of regular triangulations.
\end{proof}

In light of the preceding lemma, the minimal counterexamples necessary
to complete the proof of Theorem \ref{seccyc} are provided in our next result.

\begin{lemma} \label{compucnd} Each of the polytopes $C(9,3)$, $C(9,4)$ and
$C(9,5) $ has a triangulation and two suitable choices of points 
along the moment curve that make the triangulation regular and non-regular, 
respectively.
\end{lemma}

\begin{proof} We exhibit explicit triangulations of $C(9,3)$, $C(9,4)$ and 
$C(9,5)$ that are regular and non-regular, depending upon the choice of 
parameters.

In the case of $C(9,5)$ there is a triangulation which is regular or 
non-regular for the parameters $[0,6,7,8,9,10,11,12,30]$ and 
$[1,2,3,4,5,6,7,8,9]$, respectively:
\vskip 12pt

\noindent 
$$125689,126789,345679,125678,123489,124578,123478,124589,123457,123567$$
$$134567,256789,235679,234579,234789,245789.$$
\vskip 12pt
For the typical parameters $[1,2,3,4,5,6,7,8,9]$ the polytope $C(9,4)$ 
has 4 non-regular triangulations (out of 357), one of which is given by
the simplices
\vskip 12pt
\noindent $$34789,23789,12789,12345,46789,45678,45689,12356,$$
$$12379,12367,13479,13456, 13467,14679,14569.$$
\vskip 12pt
\noindent
On the other hand, for the parameters $[0,1/20,1/3,4,50,60,67,68,69]$ the same
triangulation becomes regular. Finally, in the case of $C(9,3)$, 
the triangulation 
$$
2578,1345,1256,1267,1278,4589,3489,2389,1289,2567,5789,3458,2358,5679,1235
$$
is non-regular for the standard parameters $[1,2,3,4,5,6,7,8,9]$ 
but becomes regular for the parameters $[1,2,3,10/3,23/6,13/3,14/3,5,6]$.
\end{proof}

\begin{remark} \rm \ \\
\noindent
One might ask whether there is a subdivision of $C(n,d)$ that is non-regular (i.e. a subdivision which is $\pi$-induced but not $\pi$-coherent for the
projection $\pi:C(n,n-1)\to C(n,d)$) for all choices of the
parameters $t_1 < \cdots < t_n$.  Such an example was recently
provided by J. Rambau (see \cite{RaSa}), who found
$4$ triangulations of $C(11,5)$ each having the property that
it is adjacent to only $4$ other triangulations by bistellar operations.
Since a triangulation which is regular for some choice of parameters
would have at least $n-1-d=11-1-5=5$ other regular neighboring vertices in
the secondary polytope, such a triangulation can never be
regular.  This is particularly interesting because of a recent result
\cite{RaSa} stating that all triangulations of $C(n,d)$ are {\it
lifting} triangulations (see \cite[\S 9.6]{OMbook} for a definition).
Any triangulation of $C(n,d)$ which is regular for some choice of
points on the moment curve is automatically a lifting triangulation,
but Rambau's examples show that the converse does not hold.
\end{remark}

\bigskip

We close our discussion by presenting in Table \ref{triangulation-numbers}
the numbers of triangulations of cyclic polytopes known to us.
Those marked with * have
been computed by J\"org Rambau. 

{\scriptsize
\begin{table}[hptb]
\label{triangulation-numbers}
\begin{center}
\begin{tabular}{| l| l | l | l | l | l | l | l | l | l | l | l |} \hline
 number of points:
              &3    & 4  & 5 & 6 & 7 & 8 & 9 & 10 & 11     & 12 \\ \hline \hline
 dimension 2  &1    & 2  & 5 &14 & 42&132&429&1430& 4862   & 16796 \\ \hline
 dimension 3  &     & 1  & 2 & 6 & 25&138&972&8477& 89405* & 1119280* \\ \hline
 dimension 4  &     &    & 1 & 2 & 7 &40 &357&4824& 96426  & 2800212*\\ \hline
 dimension 5  &     &    &   & 1 & 2 & 8 & 67&1233& 51676* & 5049932*\\ \hline
 dimension 6  &     &    &   &   & 1 & 2 & 9 & 102& 3278   & 340560*\\ \hline
 dimension 7  &     &    &   &   &   & 1 & 2 & 10 & 165    & 12589 \\  \hline
 dimension 8  &     &    &   &   &   &   & 1 & 2  & 11     & 244\\  \hline
 dimension 9  &     &    &   &   &   &   &   & 1  & 2      & 12\\ \hline
 dimension 10 &     &    &   &   &   &   &   &    & 1      & 2\\ \hline
\end{tabular}
\caption{The number of triangulations of $C(n,d)$ for $n \leq 12$.}
 \label{tricnd}
\end{center}
\end{table}
}

\section{The case $d, d'-d, n-d' \geq 2$}
\label{counterexamples}

So far we have proved all
the assertions of Theorem \ref{main-theorem} in the cases
$d=1$ (Theorem \ref{monopp})
and $n-d'=1$ (Theorem \ref{seccyc}).
Since the case $d'-d=1$ is trivial (see Section \ref{intro})
it only remains to deal with the case where $d, d'-d, n-d'$ are
all at least $2$.  We collect the assertions of Theorem 1.1 which
cover this case in the following result.
\begin{theorem}
\label{remaining-case}
If $d, d'-d, n-d' \geq 2$ then for the natural projection
$\pi: C(n,d') \rightarrow C(n,d)$ there exists a
$\pi$-induced polytopal subdivision of $C(n,d)$ whose $\pi$-coherence
depends upon the choice of points on the moment curve,
and for every choice of points there exists some
$\pi$-induced but not $\pi$-coherent subdivision.
\end{theorem}

Our proof of Theorem \ref{remaining-case} proceeds
in three steps:\\[.1in]

\noindent
STEP 1. We show that for $\pi:C(n,n-2) \rightarrow C(n,2)$
with $n \geq 6$, there is a particular $\pi$-induced subdivision
of $C(n,2)$ whose $\pi$-coherence depends upon the choice of parameters.\\

\noindent
STEP 2. We use the subdivision from Step 1 to produce a
subdivision with the same property for $\pi:C(n,d') \rightarrow C(n,2)$
whenever $d'-2, n-d' \geq 2$.\\

\noindent
STEP 3.  We use the subdivision from Step 2 to produce a subdivision
with the same property for $\pi:C(n,d') \rightarrow C(n,d)$
whenever $d, d'-d, n-d' \geq 2$.\\[.1in]

Before we proceed, we review some facts about the facial
structure of the cyclic polytopes $C(n,d)$.  From now on, we will
refer to a vertex $v_i = (t_i,t_i^2,\ldots,t_i^d)$ by its
index $i$, so that a subset $S \subseteq [n]:=\{1,2,\ldots,n\}$
may or may not span a boundary face of $C(n,d)$.  Gale's
famous Evenness Criterion \cite[Theorem 0.7]{Ziegler},\cite[\S 4.7]{Grunbaum}
tells us exactly when this happens.  The criterion is based on the
unique decomposition of $S=Y_1 \cup X_1 \cup X_2 \cup \cdots \cup X_t \cup Y_2$
of $S$ in which all $X_i,Y_j$ are contiguous segments of integers and
only $Y_1, Y_2$ may contain $1,n$, respectively (so that
$Y_1, Y_2$ may be empty).

\begin{theorem}(Gale's Evenness Criterion)
A subset $S \subseteq [n]$ spans a boundary $(|S|-1)$-face
of $C(n,d)$ if and only if in the above decomposition of $S$,
the number of interior components $X_i$ with odd length is
at most $d-|S|$.  
\end{theorem}

\subsection*{STEP 1}

The following Lemma achieves Step 1.

\begin{lemma}
\label{Step1-lemma}
For $n \geq 6$, the subdivision $T$ of $C(n,2)$ into the 
polygons $P_1=\{2,3,4,5\}$ and $P_2=\{1,2,5,6,7,\ldots,n-1,n\}$
is $\pi$-induced for $\pi:C(n,n-2) \rightarrow C(n,2)$ but
its $\pi$-coherence depends upon the choice of parameters.
\end{lemma}
\begin{proof}
Let $p_1,\ldots,p_n$ denote
the vertices of $C(n,n-2)$ and $q_1,\ldots,q_n$ the corresponding vertices
of $C(n,2)$. It is clear that $T$ is a polygonal
subdivision of $C(n,2)$ and one can check from Gale's Evenness
Criterion that it is $\pi$-induced from $C(n,n-2)$, if $n \geq 6$.

To prove that the $\pi$-coherence of $T$ depends upon the parameters
we use Lemma \ref{Paco-characterization}. Since $C(n,n-2)$
has only two more vertices than its dimension, there is a unique (up
to scaling) affine dependence $\sum_i c_i p_i=0$ among its vertices
whose coefficients $c_i$ are given by the formula
$$
c_i = \prod_{\scriptsize \begin{array}{c} j=1 \\ j\neq i \end{array}}^n
\frac{1}{t_j -t_i}
$$
as functions of the parameters $t_1 < \cdots < t_n$. For the purpose
of exhibiting parameters which
make $T$ either $\pi$-coherent or $\pi$-incoherent, we fix the
parameters $t_2=2, t_3=3, t_4=4, t_5=5$ and vary the rest.

If $T$ is $\pi$-coherent,
the functional $f({\bf x}) = \sum_i w_i x_i$ exhibiting its
$\pi$-coherence has the property that all the lifted points
$(q_i,w_i) \in \reals^3$ for $i$ in the polygon $P_2$ are
coplanar. We wish to show that we can furthermore assume that
$w_i=0$ for $i \in P_2$. To argue this,
note that coplanarity implies that there is an affine functional $h$ on
$\reals^2$ with the property that $h(q_i) = w_i$ for all
$i \in P_2$.  Hence the functional 
$f'=\sum_i w_i'x_i$ with $w_i'=w_i-h(q_i)$ will induce the same
subdivision $T$ and will have the property that
$w_i'=0$ for all $i \in P_2$. We further claim that $f'$ also exhibits
the $\pi$-coherence of $T$. Since $f$ does so, $f \in im(\phi_P^*)$,
where we are using the notation of Section \ref{background} with 
$P=C(n,n-2)$ and $Q=C(n,2)$. Moreover, the funcional $\sum_i h(q_i) x_i$
is the composition $h \circ \pi \circ \rho \circ \phi_P \in
im(\phi_P^*)$, where $\rho : \reals^{n-1} \to \reals^{n-2}$ simply
forgets the first coordinate. Hence
$f' \in im(\phi_P^*)$ as well and so $f'$ exhibits the $\pi$-coherence
of $T$.

Now notice that not only does $(w_1',\ldots,w_n')$ have the
property that $w'_i=0$ for all $i \in P_2$ but also, since
our choice of $t_2,t_3,t_4,t_5$ makes the quadrilateral $P_1$
into a trapezoid and since the lifted points $(q_i,w_i)$ are coplanar,
we have $w'_3=w'_4>0$. By Lemma \ref{Paco-characterization},
$T$ is $\pi$-coherent if and only if $\sum_i c_i w_i = 0$, which becomes
$c_3 + c_4 = 0$, or in other words $|c_3|=|c_4|$.

It now suffices to show that by suitable choices of the
remaining parameters $t_1,t_6,t_7,\ldots,t_n$
we can make the equation $|c_3|=|c_4|$ valid or invalid.
Let $K$ be a very large positive number and $\epsilon$ a very small
positive number. Consider the following
two situations:
\begin{enumerate}
\item[$\bullet$] Choosing $t_1=-K$ and $t_i$ close to $5+\epsilon$
for $i \geq 6$, one has that
\begin{equation*}
\begin{split}
c_3 &\text{ is approximately }\frac{1}{2 \cdot (2+\epsilon)^{n-5} (K+3)}, \\
c_4 &\text{ is approximately }-\frac{1}{2 \cdot (1+\epsilon)^{n-5} (K+4)}
\end{split}
\end{equation*}
and so $\frac{|c_4|}{|c_3|}$ is approximately $2^{n-5}$.
\item[$\bullet$] Choosing $t_1=2-\epsilon$ and $t_i$ close to $K$
for $i \geq 6$, one has that
\begin{equation*}
\begin{split}
c_3 &\text{ is approximately }\frac{1}{2 \cdot (K-3)^{n-5} (1+\epsilon)}, \\
c_4 &\text{ is approximately }-\frac{1}{2 \cdot (K-4)^{n-5} (2+\epsilon)} 
\end{split}
\end{equation*}
and so $\frac{|c_4|}{|c_3|}$ is approximately $\frac{1}{2}$.
\end{enumerate}
Since $\frac{|c_4|}{|c_3|}$ is greater than
$1$ in the first case but less than $1$ in the second case,
we conclude that for the choice of parameters in either
of these two situations, $T$ is not $\pi$-coherent.  However  if one
varies $t_1,t_6,t_7,\ldots,t_n$ continuously,
there must be some choice for which $\frac{|c_4|}{|c_3|} = 1$. This choice
makes $T$ $\pi$-coherent.
\end{proof}

\subsection*{STEP 2}

For Step 2, we wish
to classify the boundary faces $F$ of $C(n,d')$ into three types:
\textit{upper, lower} and \textit{contour faces} according to whether
the normal cone to $F$ in $(\reals^{d'})^*$ only contains functionals
$f({\bf x})=\sum_i w_i x_i$ in which $w_{d'}$ is positive or only
those with $w_{d'}$ negative or both kinds, respectively.
Equivalently, contour faces of $C(n,d')$ are those faces whose
normal cone contains a functional $f({\bf x})=\sum_i w_i x_i$
with $w_{d'}=0$, which is to say that they are the faces which
project to boundary faces of $C(n,d'-1)$. 

Observe that the surjection $\reals^{d'}\to ker(\pi)^*$ associated to
the projection $\pi: C(n,d')\to C(n,d)$ is the one which forgets the
first $d$ coordinates. In particular, the projections to $ker(\pi)^*$
of the normal cones of an upper face and a lower face of $C(n,d')$ do
not intersect and, as a consequence of Lemma
\ref{coherence-by-normal-cones}, any $\pi$-induced subdivision of
$C(n,d)$ having a lower and an upper face of $C(n,d')$ has to be
$\pi$-incoherent.

One can explicitly say which \textit{facets} (maximal faces) of
$C(n,d')$ are lower and which are upper
(they are never contour faces). Gale's Evenness
Criterion tells us that a subset $S \subset [n]$ forms
a facet of $C(n,d)$
if and only if all of its internal contiguous segments $X_i$ are
of even length in the decomposition 
$S=Y_1 \cup X_1 \cup X_2 \cup \cdots \cup X_t \cup Y_2$. One can then check
that the upper (respectively lower) facets are those in which
$Y_2$ has odd (respectively even) length. A non-maximal boundary face $F$
is then upper (respectively lower) if and only if it lies only in upper 
(respectively lower) facets. Otherwise $F$ is a contour face.

The point of introducing this terminology is the following
Lemma.

\begin{lemma} \label{upper}
Let $T$ be a $\pi$-induced polytopal subdivision for
$\pi:C(n,d') \rightarrow C(n,d)$ and some choice of
parameters $t_1 < \cdots < t_n$.  Assume $T$ contains
no face $\pi(F)$, where $F$ is an upper face of $C(n,d')$.
Let $T'$ be the  subdivision of $C(n+1,d)$ with parameters
$t_1 < \cdots < t_n < t_{n+1}$ obtained by adding to the
faces of $T$ all simplices obtained by adjoining the vertex
$n+1$ to the upper facets of $C(n,d)$.  Then
\begin{enumerate}
\item[(i)] $T'$ is $\pi$-induced and contains no face $\pi(F)$
for any upper faces $F$ of $C(n+1,d')$.
\item[(ii)] If $T$ is $\pi$-coherent and induced by
a functional $f({\bf x})=\sum_{i=1}^{d'} w_i x_i$ with $w_{d'} <0$
then $T'$ is $\pi$-coherent and induced by the same $f$
for any sufficiently large choice of the parameter $t_{n+1}$.
\item[(iii)] If $T$ is not $\pi$-coherent then neither is
$T'$.
\end{enumerate}
\end{lemma}
\begin{proof}
(i) This follows from the fact that every lower facet of
$C(n,d')$ is a lower facet of $C(n+1,d')$ as well and that
when one adjoins the vertex $n+1$ to an upper facet of $C(n,d)$
one obtains a (lower) facet of $C(n+1,d+1)$ and, thus, a contour face
of $C(n+1,d')$ for every $d'\ge d+2$.

(ii) According to Lemma \ref{coherence-by-normal-cones},
saying that $T$ is $\pi$-coherent and induced by
$f({\bf x})=\sum_{i=1}^{d'} w_i x_i$ is equivalent to saying
that for every cell $F$ in $T$, the functional 
$f(0,\ldots,0,x_{d+1},x_{d+2},\ldots,x_{d'})$ is in the projection
under $(\reals^{d'})^* \twoheadrightarrow ker(\pi)^*$
of the normal cone to the face of $C(n,d')$
corresponding to $F$.

If $t_{n+1}$ is chosen large, the lower faces of
$C(n+1,d')$ that use the point $n+1$, in particular the lower faces of
$C(n+1,d')$ that occur in $T'$ but were not already in $T$, are
``almost vertical''. As a consequence, the vectors with $d'$-th
coordinate negative lying in the normal cones to these new lower
facets are ``almost horizontal''. Since $w_{d'} <0$, given any cell in
$T$, we can choose $t_{n+1}$ so large that 
$f(0,\ldots,0,x_{d+1},x_{d+2},\ldots,x_{d'})$ still lies in the
projection of the normal cone of the face of $C(n+1,d')$ which corresponds to 
that cell. Thus, the coherent subdivision of $C(n+1,d)$ induced by
$f$ will contain all the cells of $T$, so it must be precisely $T'$.

(iii) Adding the extra point $n+1$ makes the normal cones of the faces
of $C(n,d')$ corresponding to cells of $T$ smaller. Since $T'$ 
contains all the cells of $T$, Lemma \ref{coherence-by-normal-cones}
implies that if $T$ were $\pi$-incoherent then $T'$ would also be
$\pi$-incoherent.  
\end{proof}

The next corollary is immediate from the previous Lemma.

\begin{corollary}
Let $T$ be a subdivision of $C(n,2)$ which is $\pi$-induced
for $\pi:C(n,d') \rightarrow C(n,2)$ with $d'\ge 4$.
Assume $T$ has some cell which corresponds to a lower face
and no cell which corresponds to an upper face of $C(n,d')$. Let $T'$ be the 
subdivision of $C(n+1,d')$ obtained adding the triangle $\{1,n,n+1\}$.
Then $T'$ satisfies the above hypotheses with $n$ replaced
by $n+1$ and
\begin{enumerate}
\item[(i)] If the $\pi$-coherence of $T$ depends upon the
choice of parameters, the same is true for $T'$.
\item[(ii)] If $T$ is $\pi$-incoherent for every choice of parameters,
then so is $T'$.
\end{enumerate}
\end{corollary}

\begin{proof}
Lower faces of $C(n,d')$ are lower faces of $C(n+1,d')$ as well
and $\{1,n,n+1\}$ is a contour face of $C(n+1,d')$ for $d'\ge
4$. Thus, $T'$ satisfies the hypotheses. Parts (i) and (ii) follow trivially
>from Lemma \ref{upper}. Part (ii) of the lemma applies to our case
because the fact that $T$ has a lower face of $C(n,d')$ implies that
\textit{any} functional $\sum_{i=1}^{d'}w_ix_i$ which
exhibits the  $\pi$-coherence of $T$ has $w_{d'}<0$. 
\end{proof}

We can now complete Step 2.  It is easy to check that the subdivision
$T$ of $C(n,2)$ produced in Lemma \ref{Step1-lemma} satisfies the
hypotheses of the previous corollary with $d'=n-2$, if $n \geq 6$:
the polygon $P_2=\{1,2,5,6,7,\ldots,n-1,n\}$ corresponds
to a lower facet of $C(n,n-2)$ and  the
polygon $P_1 = \{2,3,4,5\}$ corresponds to a lower facet of $C(6,4)$
and to a contour face of $C(n,n-2)$ for $n \geq 7$.
Therefore by iterating part (i) of the corollary we obtain
for each $d' \geq 4, n-d' \geq 2$ a subdivision $T$ of $C(n,2)$ which is 
$\pi$-induced for $\pi: C(n,d') \rightarrow C(n,2)$, but whose
$\pi$-coherence depends upon the choice of parameters.

\subsection*{STEP 3}

Here we make use of 
Lemma \ref{fiber-lifting-lemma}.
As was said in Section \ref{background}, the natural projections
$\pi:C(n,d') \rightarrow C(n,d)$ and 
$\hat{\pi}:C(n+1,d'+1) \rightarrow C(n+1,d+1)$
are compatible with the single element liftings $C(n+1,d'+1), C(n+1,d+1)$
of $C(n,d'), C(n,d)$, respectively. Recall from that section that
this required us to choose the parameters $t_1 < \cdots < t_n < t_{n+1}$
so that $t_{n+1}=0$. Again this is not a problem, as an affine 
transformation $t_i \mapsto at_i+b$ leads to compatible
affine transformations of $C(n,d),C(n+1,d+1),C(n,d'),C(n+1,d'+1)$
and $\pi$-coherence of subdivisions is easily seen to be preserved
by such transformations. Since this situation makes
$C(n,d)$ the \textit{vertex figure} \cite[\S 2.1]{Ziegler}
of $C(n+1,d+1)$ for the vertex $n+1$, any polytopal subdivision $T'$
of $C(n+1,d+1)$ gives rise to a polytopal subdivision $T$ of $C(n,d)$
by taking the \textit{link} of $n+1$ in $T'$. In this situation
we say that $T'$ \textit{extends} $T$.

\begin{proposition}
\label{step3}
Let $T$ be a polytopal subdivision of $C(n,d)$ which is
$\pi$-induced for the projection $\pi:C(n,d') \rightarrow C(n,d)$
and some choice of parameters $t_1 < \cdots < t_n$.
\begin{enumerate}
\item[(a)] If $T$ is $\pi$-coherent then 
for every choice of the parameter $t_{n+1}$, $T$ extends to a
subdivision $T'$ of $C(n+1,d+1)$ which is $\hat{\pi}$-coherent.
\item[(b)] If $T$ is not $\pi$-coherent then it does not extend to any
subdivision $T'$ of $C(n+1,d+1)$ which is $\hat{\pi}$-coherent.
\end{enumerate}
\end{proposition}
\begin{proof}
If $T$ is $\pi$-coherent then by Theorem \ref{subspace-coherence}
there is a functional 
$f \in  im(\phi^*_{P,Q}) \subset ker(\phi_Q)^*$ which
induces it. Under the identification given by the isomorphism in
Lemma \ref{fiber-lifting-lemma}, this same functional $f$ will induce
some $\hat{\pi}$-coherent subdivision $T'$ that extends $T$.
Conversely, if $T'$ were $\hat{\pi}$-coherent and extended $T$, then
the vector 
$f' \in im(\phi^*_{\hat{P},\hat{Q}}) \subset ker(\phi_{\hat{Q}})^*$
which induces $T'$ would map under the reverse of this isomorphism 
to a vector that induces $T$ and would demonstrate its $\pi$-coherence. 
\end{proof}

This finally gives us the result needed to complete Step 3.

\begin{corollary}
If $\pi: C(n,d') \rightarrow C(n,d)$ has a $\pi$-induced
polytopal subdivision of $C(n,d)$ whose $\pi$-coherence depends
upon the choice of parameters then so does
$\hat{\pi}: C(n+1,d'+1) \rightarrow C(n+1,d+1)$. 
\end{corollary}
\begin{proof}
Given such a subdivision $T$ of $C(n,d)$ and a choice of
parameters which makes it $\pi$-coherent, part (a) of the
previous corollary produces a subdivision $T'$
of $C(n+1,d+1)$ which is $\hat{\pi}$-coherent for some
choices of the parameters.  But if one chooses the first $n$
parameters $t_1 < \cdots < t_n$ so that $T$ is $\pi$- incoherent
then $T'$ must also be $\hat{\pi}$-incoherent by part (b),
regardless of how $t_{n+1}$ is chosen.
\end{proof}

\section{An instance of the Baues problem}
\label{Baues-section}

The discussion in Step 1 of Section \ref{counterexamples}
shows that the case of the projections $\pi:C(n,n-2) \rightarrow C(n,2)$
plays a special role in the family of projections 
$\pi:C(n,d') \rightarrow C(n,d)$, in that they provide examples
of non-canonical behavior ($\pi$-induced subdivisions whose 
$\pi$-coherence depends upon the choice of parameters) with
$d$ and $n-d'$ both minimal. This prompted our study of the
Generalized Baues Problem (or GBP) in this case, leading 
to Theorem \ref{new-Baues-case}. Let $\bar{\Baues}(P \overset{\pi}
{\rightarrow} Q)$ denote the proper part of the Baues poset $\Baues(P
\overset{\pi}{\rightarrow} Q)$, i.e. $\Baues(P\overset{\pi}{\rightarrow} Q)$ 
with its maximal element removed. Recall that the GBP asks whether 
$\bar{\Baues}(P \overset{\pi}{\rightarrow} Q)$ has the
homotopy type of a $(\dim(P) - \dim(Q)-1)$-sphere. When we refer
to the topology of a poset $\LL$ we always mean the topology
of the \textit{geometric realization} $|\Delta(\LL)|$
of its \textit{order complex} $\Delta(\LL)$, 
so that $\Delta(\LL)$ is the simplicial complex of chains in
$\LL$ (see \cite{Bjorner}). We now restate and prove Theorem 
\ref{new-Baues-case}.

\vspace{0.1 in}
\noindent
{\bf Theorem \ref{new-Baues-case}.} 
{\it
Let $\pi: P \rightarrow Q$ be a linear surjection of polytopes with
the following two properties:
\begin{enumerate}
\item[$\bullet$] $P$ has $\dim(P) + 2$ vertices and
\item[$\bullet$] the point configuration $\AAA$ which
is the image of the set of vertices of $P$ under $\pi$ has
only coherent subdivisions.
\end{enumerate}
Then the GBP has a positive answer for $\pi: P \rightarrow Q$,
i.e. the poset $\bar{\Baues}(P\overset{\pi}{\rightarrow} Q)$ of all
proper $\pi$-induced subdivisions of $Q$ has the homotopy type of
a $(\dimens(P)-\dimens(Q)-1)$-sphere.
}
\begin{proof}
Consider the tower of projections
$$
\Delta^{n-1} \rightarrow P \overset{\pi}{\rightarrow} Q,
$$
where $\Delta^{n-1}$ is the standard $(n-1)$-simplex in $\reals^n$,
as in Section \ref{background}. There is an obvious
inclusion
$$
\Baues(P \overset{\pi}{\rightarrow} Q) \hookrightarrow
\Baues(\Delta^{n-1} \rightarrow Q) 
$$
which simply identifies every $\pi$-induced subdivision of $Q$ with
a subdivision of $Q$ which uses only the points of $\AAA$ as
vertices. Since all subdivisions of the point configuration $\AAA$ are
coherent, $\bar{\Baues}(\Delta_{n-1} \rightarrow Q)$ can be described as
follows.  Let $Q^*=\{q_1^*,\ldots,q_n^*\}$ be the Gale transform of
$Q$. Every subdivision of $Q$ corresponds uniquely to a face
of the {\it chamber complex} of $Q^*$, i.e. to a non-zero cone in the
common refinement of all the simplicial fans generated by vectors of
$Q^*$. Thus, the poset $\bar{\Baues}(\Delta_{n-1} \rightarrow Q)$ is the
opposite or dual poset $\LL^{opp}$ to the poset $\LL$ of faces
in the  chamber complex of $Q^*$. Equivalently,
$\bar{\Baues}(\Delta_{n-1} \rightarrow Q)$ is isomorphic to the poset
of proper faces of the \textit{secondary polytope} of $Q$ (see \cite[Chapter
7]{GKZbook} or \cite{BFS} for details).

Only certain of the cones in the chamber complex of $Q^*$
correspond to $\pi$-induced subdivisions of $Q$, that is, to
elements of the subposet $\LL':=\bar{\Baues}(P \overset{\pi}{\rightarrow}
Q)^{opp} \subseteq \LL$.  We know that $\LL$, considered as a
topological space, is homeomorphic to a $(\dim(P)-\dim(Q))$-sphere and
wish to show that the subspace $\LL'$ is homotopy equivalent to a
$(\dim(P)-\dim(Q)-1)$-sphere. We will show that the cones
corresponding to the elements of $\LL-\LL'$
form two disjoint, convex (but not necessarily closed) unions
of cones $U^+, U^-$ and that there is a linear functional $f$ which separates
them, in the sense that $f(x) > 0$ for all $x \in U^+$ and $f(x)<0$ for
all $x \in U^-$. The result then follows immediately from
the technical Lemma \ref{cone-topology-lemma}.

Since $P$ has $dim(P)+2$ vertices, its vertices contain a unique
(up to a scaling factor) affine dependence
\begin{equation}
\label{unique-affine-dependence}
\sum_{i \in F^+} c_i p_i = \sum_{j \in F^-} c_j p_j
\end{equation}
with $c_i,c_j > 0$ for all $i,j$. Let $F^0=\{1,\dots,n\}\setminus
(F^+\cup F^-)$. Observe that the affine dependence projects to an
affine dependence in $Q$ and induces a functional $f$ in the Gale
transform $Q^*$ such that $f(q^*_i)$ is zero, positive or negative if
$i$ is in $F^0$, $F^+$ or $F^-$, respectively. In oriented matroid
terms, the affine dependence is a \textit{vector} in $Q$ and, thus, a
\textit{covector} in $Q^*$. This $f$ will be the linear functional of 
the statement of Lemma \ref{cone-topology-lemma}.

A subset $F$ of indices represents a face of $P$ if and only if it
either contains $F^+\cup F^-$ or contains neither $F^+$ nor $F^-$.
Thus, $F$ represents a non-face if it contains $F^+$ or $F^-$, but not
both. The complements of these non-faces are the ``forbidden cones''
in the Gale transform, according to Lemma \ref{OM-duality-lemma}.
Hence, the forbidden cones are of the following two types:
\begin{enumerate}
\item[$(+)$] $A^+ \cup A^0$ where $\emptyset\neq A^+ \subseteq F^+$,
$A^0\subseteq F^0$, or
\item[$(-)$] $A^- \cup A^0$ where $\emptyset\neq A^- \subseteq F^-$,
$A^0\subseteq F^0$.
\end{enumerate}

We claim that the unions of the forbidden open cones
of the types $(+)$ and $(-)$ are:
$$
U^+=pos(F^0\cup F^+)\cap \{f(x)> 0\} \quad \textrm{and} \quad
U^-=pos(F^0\cup F^-)\cap \{f(x)< 0\}, 
$$
respectively. Indeed, let $v$ be a vector in the relative interior of
a cone of type $(+)$. Clearly, $v$ lies in $pos(F^0\cup F^+)$ and $f(v)>0$
since every cone of type $(+)$ contains among its generators a vector
on which $f$ is positive and no vector on which $f$ is
negative. Conversely, if $v\in U^+\subset pos(F^+\cup F^0)$ then
$v$ lies in a certain cone $C$ generated by a subset of $F^0\cup F^+$
and this subset must contain an element of $F^+$ or otherwise $f(v)=0$.
A similar argument applies to $U^-$.

The sets $U^+$ and $U^-$, as defined above, are clearly convex
and separated by the functional $f$, so the proof is complete.
\end{proof}

The following technical lemma was needed in the preceding proof.
\begin{lemma}
\label{cone-topology-lemma}
Let $\F$ be a complete polyhedral fan in $\reals^d$,
i.e. a collection $\{C\}$ of relatively open
polyhedral pointed cones covering $\reals^d$, each having
its vertex at the origin, such that $C \cap C'$ is a boundary
face of $C,C'$ for each pair of cones.

Let $\LL$ be the poset of non-zero cones in $\F$ and
$\LL'$ the subposet corresponding to the elements in 
two non-empty convex (but not necessarily closed) unions
of cones $U^+, U^-$.  Assume that the cones $U^+,U^-$
are separated by some functional $f \in \reals^d$, in the sense
that $f(x)>0$ for $x \in U^+$ and $f(x)<0$ for $x \in U^-$.

Then $\LL-\LL'$ is homotopy equivalent to a $(d-2)$-sphere. 
\end{lemma}
\begin{proof}
Note that $\LL$ is the face poset of the regular (actually polyhedral)
cell complex obtained by intersecting the cones in the fan $\F$
with the unit sphere $\SU^{d-1}$. The subposet $\LL-\LL'$
then indexes the cells in the cones which cover the complement
$X:= \SU^{d-1} - (U^+ \cup U^-)$. If we show that
$X$ is homotopy equivalent to $\SU^{d-2}$ then 
the proof follows immediately by applying  
Lemma \ref{cell-complex-subspace} below with $\LL''=\LL-\LL'$.

Let $\SU^{d-2}$ be the equatorial sphere
$\SU^{d-1} \cap \{x:f(x)= 0\}$ defined by the function $f$ and $H^+$ be
the ``upper" hemisphere $H^+:=\SU^{d-1} \cap \{x:f(x) \geq 0\}$
(similarly define $H^-$).  It suffices to show that
$\SU^{d-2}$ is a deformation retract of $H^+ - U^+$
and similarly for $H^- - U^-$,
since then we can retract $X$ onto $\SU^{d-2}$ by first
retracting $H^+ - U^+$ onto $\SU^{d-2}$, keeping $H^+$ fixed,
and then retracting $H^- - U^-$ onto $\SU^{d-2}$.
To retract $H^+ - U^+$ onto $\SU^{d-2}$, note that the pair
$(H^+,U^+)$ is homeomorphic to the pair $(\B^{d-1},C)$ by forgetting
the last coordinate and then scaling, where $C$ is some convex subset
inside the unit disk $\B^{d-1}$. We then need to retract $\B^{d-1}-C$
onto the boundary $\partial \B^{d-1} = \SU^{d-2}$,
which can be done as follows. Pick any point $p \in C$.
For $x \in \B^{d-1}-C$, let $s(x)$ be the unique
point on the boundary $\partial \B^{d-1}$ that lies on the ray
which emanates from $p$ and passes through $x$. Define the
deforming homotopy
$f:(\B^{d-1}-C) \times [0,1] \rightarrow \B^{d-1}-C$
by $f(x,t) = (1-t) \, x + t \, s(x)$.  Convexity of $C$ shows that $f$ is
well-defined, i.e. that its image lies in $\B^{d-1}-C$.  Also $f$
is continuous because $s$ is the restriction to $\B^{d-1}-C$
of a continuous map $\B^{d-1}-\{p\} \rightarrow \partial\B^{d-1}$.
\end{proof}

The following lemma is probably well-known but we do not
know of a proof in the literature, so we include one here.

\begin{lemma}
\label{cell-complex-subspace}
Let $K$ be a regular cell complex with face poset $\LL$ and let $\sigma_l$
be the cell of $K$ indexed by $l \in \LL$. Then for any subposet 
$\LL'' \subseteq \LL$,
the subspace $K'' : = \cup_{l \in \LL''} \sigma_l$ is
homotopy equivalent to $\LL''$.
\end{lemma}
\noindent
\begin{proof}
Since $K$ is a regular cell complex and $\LL$ is its face poset,
the order complex $\Delta(\LL)$ is the first barycentric subdivision
$Sd(K)$ and the order complex $\Delta(\LL'')$ is a subcomplex
of $Sd(K)$.  The subspace $K''$ of $K$ may be identified with
a subspace of $Sd(K)$. We will describe a deformation retraction
of $K''$ inside $Sd(K)$ onto the subcomplex $\Delta(\LL'')$.
The retraction will be defined piecewise on each simplex of $Sd(K)$.
A simplex $\sigma$ of $Sd(K)$ is represented by a chain
$l_1 < \cdots < l_r$ in $\LL$ and has vertices labelled $l_1,\ldots,l_r$.
Let $\sigma_1$ be the subface of this simplex spanned by the $l_i$'s
which lie in the subposet $\LL''$ and $\sigma_2$ be the opposite subface, i.e.
the one spanned by the rest of the $l_i$. Then $|\sigma|$ is
the \textit{topological join} 
$$
|\sigma| = |\sigma_1| * |\sigma_2| 
: = |\sigma_1| \times |\sigma_2| \times [0,1] / ((x,y,0) \sim x, (x,y,1) \sim y).
$$
One can check from the definition that
$|\sigma \cap K''| \subseteq |\sigma|-|\sigma_2|$
and $|\sigma \cap \Delta(\LL'')| = |\sigma_1|$.
Since for any topological join $X * Y$ one can retract
$X * Y - Y$ onto $X$, we can retract $|\sigma \cap K''|$ onto
$|\sigma \cap \Delta(\LL'')|$ for each simplex $\sigma$.  It is easy
to see that all of these retractions can be done coherently, giving
a retraction of $K''$ onto $|\Delta(\LL'')|$, as desired.
\end{proof}

As was mentioned earlier, Theorem \ref{new-Baues-case}
has the following corollary.

\begin{corollary}
The GBP has a positive answer for $\pi: C(n,n-2) \rightarrow C(n,2)$.
\end{corollary}

\bigskip

It is perhaps worthwhile to look more closely at the first
interesting example, i.e. the projection $\pi:C(6,4) \rightarrow
C(6,2)$. From the results of Section \ref{counterexamples}, this is the
minimal example where the fiber polytope depends upon the choice of
parameters. We first compute the Baues poset $\Baues(C(6,4)
\rightarrow C(6,2))$ (which does not depend on the choice of
parameters) by using the technique in the proof of Theorem
\ref{new-Baues-case}. The secondary polytope of $C(6,2)$ is the
well-known 3-dimensional associahedron, which is a simple 3-polytope
with six pentagons and three quadrilaterals as facets. A picture of it
can be found in \cite[page 239]{GKZbook}. There are two special
vertices, each incident to three pentagons, which correspond to the
triangulations $\{135,123,345,561\}$ and $\{246,234,456,126\}$. In the
chamber complex of the Gale transform, these triangulations correspond
to chambers that are triangular cones generated by $246$ and $135$,
respectively. The other twelve vertices of the associahedron
are incident to two pentagons and a quadrilateral each.

\begin{figure}[hptb]
\center{\mbox{\epsfbox{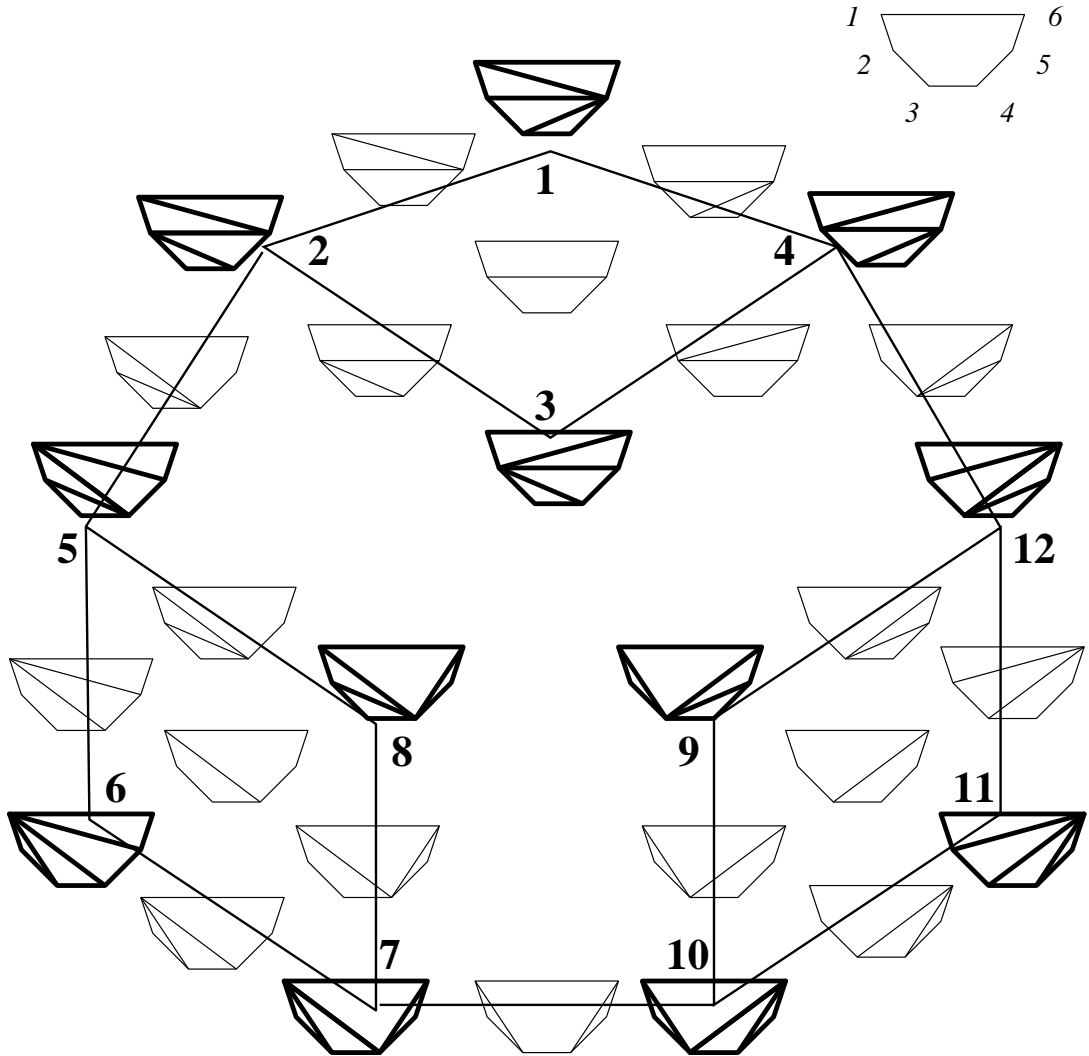}}} 
\caption{The structure of $\pi$-induced subdivisions.}
\label{pi}
\end{figure}

The minimal non-faces of $C(6,4)$ are $F^+=135$ and $F^-=246$. Since in this
case $F^0$ is empty (this is always the case in $C(n,n-2)$),
the regions $U^+$ and $U^-$ in the proof of Theorem
\ref{new-Baues-case} to be removed from the chamber complex of
$C(6,2)$ are the closed triangular cones generated by
$135$ and $246$. In other words, the
Baues poset $\Baues(C(6,4) \rightarrow C(6,2))$ is isomorphic to the
poset of proper faces of the associahedron not incident to the two
special vertices mentioned above. This leaves us with twelve vertices
of the associahedron, representing twelve $\pi$-induced
triangulations, 15 edges, representing 15 $\pi$-induced subdivisions
of height one in the poset and 3 quadrilaterals, representing 3
$\pi$-induced subdivisions of height 2.
The cell complex of these faces is depicted in Figure
\ref{pi}, where we have drawn the subdivision corresponding to each
face. The twelve subdivisions in thick are the triangulations, which
we have numbered from {\bf 1} to {\bf 12}. In the following discussion
we will refer to a $\pi$-induced subdivision by the $\pi$-induced
triangulations which refine it. Thus, {\bf (1,2,3,4)}, {\bf
(5,6,7,8)} and {\bf(9,10,11,12)} represent the three subdivisions of
height two (the quadrilaterals in Figure \ref{pi}).

We now want to study the fiber polytope associated to the projection
$\pi$ and its dependence with the choice of parameters. We first
recall that a $\pi$-coherent subdivision cannot have both lower and
upper faces of $C(6,4)$. The upper and lower facets of $C(6,4)$ are
$\{1234,1245,$ $1256,$ $2345,$ $2356,3456\}$ and
$\{1236,1346,1456\}$, respectively. Thus, the faces $136$ and $146$ are
upper, while $234$ and $345$ are lower. This implies that the following
subdivisions are not $\pi$-coherent:
$$
{\bf (8)}=\{124,234,146,456\},\qquad\qquad
{\bf (9)}=\{123,136,345,356\},
$$
$$
{\bf (7,8)}=\{1234,146,456\},\qquad
{\bf (5,8)}=\{124,234,1456\},\qquad\quad\quad
{\bf (5,6,7,8)}=\{1234,1456\},
$$
$$
{\bf (9,10)}=\{1236,345,356\},\qquad
{\bf (9,12)}=\{123,136,3456\},\qquad
{\bf (9,10,11,12)}=\{1236,3456\}.
$$

Incidentally, the same argument shows that for the projection $\pi:
C(2d'-2,d')\to C(2d'-2,2)$, $d'\ge 4$, there are $\pi$-induced
subdivisions which are $\pi$-incoherent in every choice of parameters. Namely,
the subdivision with only two cells $\{1,\dots,d'\}$ and
$\{1,d',d'+1,\dots,2d'-2\}$.  Other cases can be obtained from this by
applying Proposition \ref{step3}.  We do not know whether
$n-d', d'-d, d \geq 2$ always implies there exists a subdivision of $C(n,d)$
which is $\pi$-induced for $\pi: C(n,d')\to C(n,d)$ but
$\pi$-incoherent for every choice of parameters.

After removing those non-coherent subdivisions, the poset is already
almost the face poset
of a polygon, except for the quadrilateral {\bf (1,2,3,4)}.
In particular, the triangulations {\bf 2}, {\bf 5}, {\bf 6}, {\bf 7},
{\bf 10}, {\bf 11}, {\bf 12} and {\bf 4}
are always $\pi$-coherent as well as the height one
subdivisions (bistellar flips) joining them. The possible subdivisions
whose $\pi$-coherence may depend upon the choice of parameters are listed
in the following rows:
$$
{\bf (1,2)}=\{125,156,2345\}, \qquad
{\bf (1)}=\{125,156,235,345\}, \qquad
{\bf (1,4)}=\{1256,235,345\},
$$
$$
{\bf (1,2,3,4)}=\{1256,2345\},
$$
$$
{\bf (2,3)}=\{1256,234,245\}, \qquad
{\bf (3)}=\{126,256,234,245\}, \qquad
{\bf (3,4)}=\{126,256,2345\}.
$$

The $\pi$-coherence of the subdivision $\{1256,2345\}$ is precisely what was
studied in Step 1 of Section \ref{counterexamples} and, actually, the
three cases of $|c_3|$ being less, equal or greater than $|c_4|$
which appeared there produce the three cases for the fiber
polytope. In the two extreme cases the fiber polytope is a 9-gon and in
the middle one it is an octagon.

\section{Acknowledgments}
The authors thank Donald Kahn and Ilia Itenberg for helpful 
conversations regarding the proof of Lemma \ref{cone-topology-lemma}.
They also offer special thanks to J\"org Rambau for his amazing software
which calculated many of the entries in Table \ref{triangulation-numbers}.

\end{document}